\documentclass[11pt]{article}

\usepackage{amssymb}

\newtheorem{theorem}{$~~~~$ Theorem}[section]
\newtheorem{corollary}{Corollary}

\newtheorem{example}[theorem]{$~~~~$ Example}
\newtheorem{lemma}[theorem]{$~~~~$ Lemma}
\newtheorem{remark}[theorem]{$~~~~$Remark}

\newtheorem{definition}[theorem]{$~~~~$Definition}

\def\nor1{Normed$\{~2^{ \zzz \theta  \, )} ~$,$~\sqrt{~2^{ \zzz \theta  \, )}}~\}$}

\def\xor2{Normed$\{ ~\sqrt{~2^{ \zzz \theta  \, )}}~,~2~ \} $}

\def\zzz{~ \sharp ( ~ }

\def\f55{ \normalsize  \baselineskip = 1.8 \normalbaselineskip }

\def\f55{  \baselineskip = 1.1 \normalbaselineskip } 
\def\g55{  \baselineskip = 1.0 \normalbaselineskip } 
\def\s55{ \baselineskip = 1.0 \normalbaselineskip }

\newcommand{\thx}[1]{Theorem \ref{#1}}

\newcommand{\eq}[1]{(\ref{#1})}
\newcommand{\ep}[1]{Equation (\ref{#1})}

\newcommand{\txl}[1]{Tab$- {#1}$}

\begin{document}

 \title{On the Significance of 
Self-Justifying Axiom Systems
from the Perspective of Analytic Tableaux}

\def\aaa{A}
\def\ccc{Class}

\def\ulxyz{\lceil}
\def\urxyz{\rceil}

\def\ulxyz{\ulcorner}
\def\urxyz{\urcorner}

\def\beq{\begin{equation}}
\def\enq{\end{equation}}

\def\bel{\begin{lemma}}
\def\enl{\end{lemma}}

\def\bec{\begin{corollary}}
\def\enc{\end{corollary}}

\def\bed{\begin{description}}
\def\ennd{\end{description}}
\def\bee{\begin{enumerate}}
\def\ene{\end{enumerate}}

\def\bxbxd{\begin{definition}}
\def\bxbxdd{\begin{definition}}
\def\eedd{\end{definition}}
\def\bxbxdr{\begin{definition} \rm}
\def\bel{\begin{lemma}}
\def\enl{\end{lemma}}
\def\ent{\end{theorem}}

\author{  Dan E.Willard\thanks{This research 
was partially supported
by the NSF Grant CCR  0956495.
Email = dew@cs.albany.edu.}}

\date{State University of New York at Albany}

\maketitle

 \setcounter{page}{0}
 \thispagestyle{empty}

\normalsize

\baselineskip = 1.5\normalbaselineskip

\normalsize

\baselineskip = 1.0 \normalbaselineskip 
\def\bbint{\large \baselineskip = 1.6 \normalbaselineskip } 
\def\bbint{\large \baselineskip = 1.6 \normalbaselineskip }
\def\bbint{\normalsize \baselineskip = 1.3 \normalbaselineskip }

\def\bbint{\normalsize \baselineskip = 1.27 \normalbaselineskip }

\def\bbint{\large \baselineskip = 2.0 \normalbaselineskip }

\def\bbint{\normalsize \baselineskip = 1.25 \normalbaselineskip }
\def\bbina{\normalsize \baselineskip = 1.24 \normalbaselineskip }

\def\bbint{\large \baselineskip = 2.0 \normalbaselineskip }

\def\bbing{ }
\def\bbins{ }
\def\bbinm{ }

\def\bbint{\normalsize \baselineskip = 1.95 \normalbaselineskip }

\def\bbing{ }
\def\bbins{ }
\def\bbinm{ }

\def\bbint{\large \baselineskip = 2.3 \normalbaselineskip } 
\def\bbing{ }
\def\bbins{ }
\def\bbinm{ }

\def\bbint{\normalsize \baselineskip = 1.7 \normalbaselineskip } 

\def\bbint{\large \baselineskip = 2.3 \normalbaselineskip } 
\def\bbinm{ \baselineskip = 1.18 \normalbaselineskip }

\def\bbint{\large \baselineskip = 2.0 \normalbaselineskip } 
\def\bbing{ }
\def\bbins{ }
\def\bbinm{ }
\def\bbinr{ }

\def\bbint{\normalsize \baselineskip = 1.25 \normalbaselineskip }
\def\bbina{\normalsize \baselineskip = 1.24 \normalbaselineskip }
\def\bbinr{ \baselineskip = 1.3 \normalbaselineskip }
\def\bbing{ \baselineskip = 1.28 \normalbaselineskip }
\def\bbins{ \baselineskip = 1.21 \normalbaselineskip }
\def\bbinm{  }

\def\fgf {\baselineskip = 1.3 \normalbaselineskip }

\bbint

\parskip 5 pt

\noindent

\begin{abstract}
\footnotesize
\baselineskip = 0.9\normalbaselineskip 

This article will be a continuation of our 
research into self-justifying 
systems.
It will introduce 
several 
new theorems 
and  
their 
applications.
(One of these results 
will transform our previous infinite-sized
self-verifying 
formalisms into 
tighter
systems, 
with
only a finite number of axioms.)
It will explain how self-justification
is useful, even when the Incompleteness
Theorem 
limits its reach.

\end{abstract}

\noindent

\bigskip\bigskip

{\bf Historical Remark about this Draft:} $~$All
the theorems in this January 2014 manuscript had also
appeared in the
 {\it ``Version 1''} June 2013 manuscript
that was posted
earlier in the Cornell archives.
The 
main 
difference between these two manuscripts is that
I 
essentially
spiced up the accompanying narrative, 
in the 
 January 2014
draft,
so as to better
explain the motivation for this research endeavor. 

\bigskip\bigskip
\normalsize

\baselineskip = 1.2\normalbaselineskip \parskip 5pt 

\bigskip

{\bf Keywords:}
G\"{o}del's Second Incompleteness Theorem, Consistency, Hilbert's Second
Open Question,
Semantic Tableaux

\bigskip

{\bf Mathematics Subject Classification:}
03B52; 03F25; 03F45; 03H13 

\newpage

\def\ww22{\normalsize \baselineskip = 1.21\normalbaselineskip \parskip 4 pt}
\def\bb22{\normalsize \baselineskip = 1.19\normalbaselineskip \parskip 4 pt}
\def\zz22z{\normalsize \baselineskip = 1.19 \normalbaselineskip \parskip 3 pt}
\def\xx22{\normalsize \baselineskip = 1.17\normalbaselineskip \parskip 4 pt}
\def\vx22s{\normalsize \baselineskip = 1.16 \normalbaselineskip \parskip 3 pt} 
\def\vv22{\normalsize \baselineskip = 1.17 \normalbaselineskip \parskip 3 pt} 
\def\aa22{\normalsize \baselineskip = 1.15 \normalbaselineskip \parskip 3 pt} 
\def\f55{  \baselineskip = 1.1 \normalbaselineskip } 
\def\h55{  \baselineskip = 1.08 \normalbaselineskip } 
\def\g55{  \baselineskip = 1.0 \normalbaselineskip } 
\def\s55{ \baselineskip = 1.0 \normalbaselineskip } 
\def\sm55{ \baselineskip = 0.9 \normalbaselineskip }

\vspace*{- 1.0 em}

\def\waw11{\normalsize \baselineskip = 1.72\normalbaselineskip}
\def\waw11{\normalsize \baselineskip = 1.12\normalbaselineskip}
\def\waw11{\normalsize \baselineskip = 1.85\normalbaselineskip}

\def\waw11{\normalsize \baselineskip = 1.45\normalbaselineskip}

\def\waw11{\normalsize \baselineskip = 1.7\normalbaselineskip}

\def\waw11{\normalsize \baselineskip = 1.12\normalbaselineskip}

\def\f55{  \baselineskip = 1.59 \normalbaselineskip } 
\def\g55{  \baselineskip = 1.50 \normalbaselineskip } 
\def\s55{ \baselineskip = 1.50 \normalbaselineskip } 
\def\sm55{ \baselineskip = 1.5 \normalbaselineskip }

\def\f55{  \baselineskip = 1.59 \normalbaselineskip } 
\def\g55{  \baselineskip = 1.50 \normalbaselineskip } 
\def\s55{ \baselineskip = 1.50 \normalbaselineskip } 
\def\sm55{ \baselineskip = 0.9 \normalbaselineskip }

\def\aa22{\normalsize  \waw11 \parskip 6 pt} 
\def\bb22{\normalsize  \waw11 \parskip 5 pt}
\def\ww22{\normalsize \waw11 \parskip 4 pt}
\def\vv22{\normalsize  \waw11 \parskip 3 pt} 
\def\tt22{\normalsize  \waw11 \parskip 2 pt} 

\def\f55{  \baselineskip = 1.1 \normalbaselineskip } 
\def\g55{  \baselineskip = 1.0 \normalbaselineskip } 
\def\b55{  \baselineskip = 1.0 \normalbaselineskip } 
\def\s55{ \baselineskip = 1.0 \normalbaselineskip } 
\def\sm55{ \baselineskip = 0.9 \normalbaselineskip }

\def\mal{ \bf  }
\def\nal{\mathcal}
\def\cvt{ \baselineskip = 0.98 \normalbaselineskip }
\def\cv9{ \baselineskip = 0.99 \normalbaselineskip }
\def\cvs{ \baselineskip = 1.0 \normalbaselineskip }
\def\cvl{ \baselineskip = 1.0 \normalbaselineskip }
\def\cvh{ \baselineskip = 1.03 \normalbaselineskip }
\def\cvg{ \baselineskip = 1.00 \normalbaselineskip }

\def\cvt{ \baselineskip = 1.6 \normalbaselineskip }
\def\cv9{ \baselineskip = 1.6 \normalbaselineskip }
\def\cvs{ \baselineskip = 1.6 \normalbaselineskip }
\def\cvl{ \baselineskip = 1.6 \normalbaselineskip }
\def\cvh{ \baselineskip = 1.6 \normalbaselineskip }
\def\cvg{ \baselineskip = 1.6 \normalbaselineskip }
\def\cvb{ \baselineskip = 1.6 \normalbaselineskip }
\def\cvnew{ \baselineskip = 1.6 \normalbaselineskip }
\def\cvmew{ \baselineskip = 1.6 \normalbaselineskip }
\def\cvwew{ \baselineskip = 1.6 \normalbaselineskip \parskip 5pt }
\def\cvrew{ \baselineskip = 1.6 \normalbaselineskip \parskip 3pt }

\def\cvt{ \baselineskip = 1.22 \normalbaselineskip }
\def\cv9{ \baselineskip = 1.22 \normalbaselineskip }
\def\cvs{ \baselineskip = 1.22 \normalbaselineskip }
\def\cvl{ \baselineskip = 1.22 \normalbaselineskip }
\def\cvh{ \baselineskip = 1.22 \normalbaselineskip }
\def\cvg{ \baselineskip = 1.22 \normalbaselineskip }
\def\cvb{ \baselineskip = 1.22 \normalbaselineskip }
\def\cvnew{ \baselineskip = 1.4 \normalbaselineskip }
\def\cvmew{ \baselineskip = 1.35 \normalbaselineskip }
\def\cvwew{ \baselineskip = 1.4 \normalbaselineskip \parskip 5pt }
\def\cvrew{ \baselineskip = 1.22 \normalbaselineskip \parskip 3pt }

\def\fend{ 

\medskip -------------------------------------------------------}

\def\f55{  \baselineskip = 1.1 \normalbaselineskip } 
\def\g55{  \baselineskip = 1.0 \normalbaselineskip } 
\def\s55{ \baselineskip = 1.0 \normalbaselineskip } 
\def\sm55{ \baselineskip = 1.0 \normalbaselineskip } 
\def\h55{  \baselineskip = 1.08 \normalbaselineskip } 
\def\b55{  \baselineskip = 1.1 \normalbaselineskip } 

\def\nop{ }

\cvl
\cvnew

\cvnew
\cvl

\parskip 3pt

\def\hgskip{ \medskip }

\def\njp{\newpage}
\def\njp{ }

\cvl
\cvnew

\def\nskip{\bigskip}

\cvl
\cvnew

\def\cpl{\parskip 2pt}
\def\cpz{\parskip 1pt}
\def\cpn{\parskip 0pt}

\baselineskip = 1.04\normalbaselineskip
\parskip 3pt

\def\xxitch{\switch}

\def\switch {\normalsize  \baselineskip = 2.3\normalbaselineskip 
\parskip 7pt }

\def\switch {\normalsize  \baselineskip = 1.85\normalbaselineskip \parskip 5pt }
\def\switch {\normalsize  \baselineskip = 1.95\normalbaselineskip \parskip 5pt }

\def\switch {\large  \baselineskip = 1.95\normalbaselineskip \parskip 8pt }
\def\switch {\Large  \baselineskip = 1.95\normalbaselineskip \parskip 8pt }

\def\switch {\normalsize  \baselineskip = 1.85\normalbaselineskip \parskip 5pt }

\def\switch {\normalsize  \baselineskip = 2.1\normalbaselineskip \parskip 8pt }
\def\switch {\large  \baselineskip = 2.1\normalbaselineskip \parskip 8pt }

\def\switch {\normalsize  \baselineskip = 1.37\normalbaselineskip \parskip 6pt }

\def\switch {\large  \baselineskip = 2.1\normalbaselineskip \parskip 8pt }
\def\switch {\normalsize  \baselineskip = 1.25\normalbaselineskip \parskip 5pt }
 
\def\switch {\normalsize  \baselineskip = 1.8\normalbaselineskip \parskip 8pt }
\def\switch {\normalsize  \baselineskip = 1.4\normalbaselineskip \parskip 5pt }

\def\switch {\normalsize  \baselineskip = 1.8\normalbaselineskip \parskip 8pt }

\def\cvl{ \switch}
\def\cvu{ \switch}
\def\cvh{ \switch}
\def\cvs{ \switch}
\def\cvz{ \switch}
\def\cvm{ \switch}
\def\xxitch{\switch}
\def\britch{\switch}
\def\britch{\normalsize  \baselineskip = 1.2\normalbaselineskip }

\vspace*{- 3.0 em}

\normalsize

\cvu

\switch

\normalsize

\cvl

\vspace*{- 2.2 em}

\switch

\cvh

\bigskip
\bigskip

\bigskip
\bigskip

\section{Introduction}

\label{sss1}
\label{tst1}

\switch

\cvh

G\"{o}del's Incompleteness Theorem 
is a 2-part result.
Its
first half indicates no decision
procedure can identify all the true statements of arithmetic.
Its
 Second Incompleteness Theorem specifies
sufficiently strong systems
cannot verify their own consistency.
G\"{o}del 
was  careful to insert the following caveat
into
his historic paper 
\cite{Go31}, 
indicating  
a
diluted 
form
of Hilbert's Consistency Program 
could be successful:
\begin{quote} 
\small
\baselineskip = 0.95\normalbaselineskip
\xxitch
$~*~$ : ``It must be expressly noted that
Proposition XI 
(e.g. the Second Incompleteness Theorem)
represents no contradiction of the formalistic
standpoint of Hilbert. For this standpoint
presupposes only the existence of a consistency
proof by finite means, and { there might
conceivably be finite proofs} which cannot
be stated in P (or in M or in A). ''
\end{quote}
Yourgrau has summarized,
in detail,
G\"{o}del's 
considerations
about this 
subject.
Thus 
\cite{Yo5}   indicated that 
{\it ``for several years''} after \cite{Go31}'s publication,
G\"{o}del 
{\it ``was cautious not to prejudge''}
whether 
some
unusual formalism, 
different from
Peano Arithmetic, might provide
some
type 
of proof
of its own
consistency. 
Likewise, the Stanford Encyclopedia
\cite{Stanford}
 cites 
G\"{o}del 
remarking that it was only after 
the evolution of Turing's
formalism \cite{Tu37}
that G\"{o}del viewed   the
Second Incompleteness Theorem as 
being
fundamentally
ubiquitous.

Within the framework of these particular caveats, our
papers 
\cite{ww93}--\cite{ww11}
have discussed
both generalizations
and boundary-case exceptions for the Second Incompleteness
Theorem.
This research
has had two facets because the meaning of 
partial exceptions to the Second Incompleteness Theorem
could be easily
misunderstood, if their limitations 
are
not also
carefully 
scrutinized
and recognized
in 
meticulous detail.

\def\iii{IS$_D(\aaa)$}
\def\I2{IS$^{\#}_D(\beta)$}
\def\ik3{IS$^{\#}_D(\beta_{A,i})$}

The prior research was mostly mathematical in 
character, in that
it sought to formalize several
well-defined pairs $(T,E)$, where $~T~$ was a threshold
sufficient for enacting the force of the Second Incompleteness
Theorem and $~E~$ was a 
closely related
``boundary-case exception'' that omitted some part of $T$'s formalism.
Our goal in this paper will be 
different.
It  will
seek
to explore this subject from a  more
futuristic
and  epistemological
perspective.
It will characterize
how far a
theorem prover can traverse in understanding its own consistency
before it reaches the 
inevetible
barriers 
imposed by the 
Incompleteness Theorem.
It will
 provide some partially sympathetic
interpretations of  
Hilbert \cite{Hil25}'s
1925 statement $**$, 
$~$in a context where G\"{o}del \cite{Go31}
established it,
obviously, 
needed
 serious 
amendments:
\begin{quote} 
\xxitch
$**$
``Where 
else 
would 
reliability and truth be found 
if  even mathematical thinking fails?   The definitive nature
of the infinite has become necessary,  not merely for the special
interests of individual sciences, but rather { for the
honor} of human understanding itself.''
\end{quote}
\baselineskip = 1.02\normalbaselineskip
\parskip 1pt
\xxitch

In a context where
G\"{o}del's statement $*$ was
{\it partially}
sympathetic with Hilbert's 
goals,
 {\it even 
after }
G\"{o}del 
had proven
his incompleteness 
theorem,
we will show how 
some logic systems
can support a 
positive {\it but-curtailed}
form 
of
Hilbert's objectives.
From both a theoretical and 
a
pragmatic 
engineering-style 
perspective,
our 
results
will helpfully  explain
how to transform 
\cite{ww93,ww1,ww6,wwapal,ww9,ww11}'s 
results about
infinite-sized self-justifying
systems
into compact tableaux formalisms,
containing only a 
{\it strictly
finite} number of proper axioms.

These 
results
will suggest that although
the Second Incompleteness
Theorem is
extraordinarily robust and ubiquitous 
from an idealized purist mathematical perspective,
some unusual axiom systems 
can possess a
fragmentary
 knowledge of their consistency,
when using
definitions of consistency that are diluted
{\it but
not 
immaterial.}

This research
was partially influenced by
Goldstein's biography of G\"{o}del \cite{Go5}. It suggested
that there ought to be some type of partial
philosophical compromise available
between the mathematical Platonism of G\"{o}del and the
logical positivism of Wittgenstein.

Also, the style of presentation in this article was 
influenced by a suggestion that Selmer Bringsjord had
made
several year ago \cite{Br94}. 
It was 
that our research 
should
be 
ideally
broken into two stages, with its first phase
focusing on purely mathematical
results (similar to 
\cite{ww93}--\cite{ww11}'s treatment), and 
with
its second stage exploring the
epistemological and futuristic implications of such results.
Bringsjord felt this
carefully segmentized
 approach was 
useful
so that the  mathematical novelty
of our results
not be confused with 
its philosophical 
interpretation.
(This is because the latter,
quite naturally, 
lends itself 
much 
more easily
to  a quite complicated debate
about the underlying epistemological implications
of 
our mathematical
formalism.)

In essence, the current
paper
 will bring this second stage of
our research to completion. It will explain the significance of
logics
that
formalize  a partial,
{\it although admittedly
fragmentary,}
knowledge of their own consistency.

\medskip

{\bf Format of our Presentation:} $~$
As much as is 
feasible,
this article will be written
in 
an
informal style,  
so as to make it
comprehensible
to a
broader
audience.
A reader who
 has, 
thus, 
mastered either
say Enderton's
or Mendelson's introductory logic textbooks
\cite{End,Mend},
during a 1-semester course in  logic, 
 should be able
follow its gist.

Section \ref{mini} will suggest that  our fragmentized approach
for
interpreting the goals of 
an
applications-oriented
engineering-styled subset of mathematics
be
called 
 ``Miniaturized Finitism''.
The cautious-sounding adjective
of
{\it  ``miniaturized''} 
was 
attached to 
our formalism's
 name
partly because 
it
involves
a curtailed notion of growth
(e.g. see Sections \ref{ss7-8}
and \ref{ss8-9}) and
also because the classic variants of the Second Incompleteness
Theorem, 
obviously, have
resolved 
at least
90\% of the issues
raised by Hilbert's Second Open Question.

Our discourse
will focus, essentially,  on 
the remaining, say,
10 \% of the issues
raised by Hilbert's   year-1900
penetrating
Second Open Question. 
Thus,
we seek to 
 offer a novel  interpretation of 
G\"{o}del's and Hilbert's
1925 and 1931 
statements of $*$ and $**$,
that will
explain how human beings
manage to gain the motivation and mental energy required 
for 
stimulating their
cogitations, by having 
a type of automatic and psychologically unconscious
access to an
essentially 
{\it very very}
miniaturized
but-spontaneous self-appreciation of their own
consistency.

\section{Background Setting} 

\label{ss2-3}
\label{tst2}

Throughout this article,
$~\alpha~$ 
will
denote an axiom system, 
and $~d~$ will  
denote
 a
deduction method.
An ordered pair
 $~(  \alpha  , d  )$
will 
be called  {\bf Self Justifying} when:
\begin{description}
\xxitch
  \item[  i   ] one of $ \, \alpha \,$'s  theorems
will
state that the deduction method $ \, d, \, $ applied to the
system $ \, \alpha, \, $ will 
produce a consistent set of theorems, and
\item[  ii   ]
     the axiom system $ \, \alpha  \, $ is in fact consistent.
\end{description}
For any  $\,(\alpha,d) \,$, 
it is 
easy 
to construct a second
axiom system $ \, \alpha^d \, \supseteq  \,  \alpha  \, $
 that  satisfies
Part-i of 
this definition.
For instance,  $ \, \alpha^d \, $  could
consist of all of $~\alpha \,$'s axioms plus the following added
sentence,  that we call
{\bf SelfRef$(\alpha,d)~$}:
\topsep -3pt
\begin{quote} 
\xxitch
$\bullet~~~$ 
There is no proof 
(using 
$d$'s deduction method)
of  $0=1$
from the  {\it union}
 of
the
 axiom system $\, \alpha \, $
with {\it this}
sentence  ``SelfRef$(\alpha,d) \,$'' (looking at itself).
\end{quote}
Kleene 
\cite{Kl38} 
discussed
how
to
encode
approximate
 analogs of
SelfRef$(\alpha,d)$'s
 self-referential statement.
Each of
Kleene, 
Rogers and Jeroslow 
 \cite{Kl38,Ro67,Je71}
 noted
$\alpha ^d$ 
may,
however,  be inconsistent
(despite SelfRef$(\alpha,d)$'s assertion),
thus causing 
it
to violate Part-ii of   self-justification's
definition.

This problem arises in
settings
more general than 
 G\"{o}del's
paradigm,
where $\alpha$  was an extension of Peano Arithmetic.
There 
are 
many
settings 
where the Second Incompleteness Theorem does
generalize
\cite{Ad2,AZ1,BS76,Bu86,BI95,Go31,HP91,HB39,KT74,Lo55,Pa71,Pu85,Ro67,So94,Sv7,Vi5,WP87,ww2,wwlogos,wwapal,ww7}.
Each such result formalizes a 
paradigm where 
self-justification is infeasible,
due to a diagonalization issue.
Many
logicians 
have, thus,
hesitated 
to
 employ 
a SelfRef$(\alpha,d)$
 axiom
because
$\alpha+$SelfRef$(\alpha,d)  $
is 
typically 
inconsistent
     \footnote{ \fgf\label{troub} 
     Typically, 
     $~\alpha^d~=~\alpha \,+ \,$SelfRef$(\alpha,d)$ 
will
     be inconsistent
because
     a  
     standard
      G\"{o}del-like self-referencing 
     construction
  will 
     produce a proof of $0=1$ from
     $~\alpha^d\,$, {\it even when $~\alpha~$ is 
     consistent.}}.

Our research
explored special
circumstances
\cite{ww1,ww5,ww6,wwapal}
where it is feasible to
construct self-justifying formalisms.
These paradigms involved weakening
the properties 
a system can prove about
addition and/or 
multiplication
(to avoid the preceding
difficulties).
To be more precise, let
 $Add(x,y,z)$ and    $Mult(x,y,z)$ 
denote 
two 
3-way predicates
specifying
$x+y=z$ and
$x*y=z$.
Then a
logic
will be said to
{\bf recognize}
successor, 
 addition  and multiplication
as {\bf Total Functions} iff it 
includes
1-3 as axioms.

\vspace*{- 0.4 em}
{\small
\cvl
\beq 
\label{totdefxs}
\forall x ~ \exists z ~~~Add(x,1,z)~~
\enq
\vspace*{- 1.7 em}
\beq 
\label{totdefxa}
\forall x ~\forall y~ \exists z ~~~Add(x,y,z)~~
\enq
\vspace*{- 1.7 em}
\beq 
\label{totdefxm}
\forall x ~\forall y ~\exists z ~~~Mult(x,y,z)~
\enq }

\vspace*{- 1.2 em}

\xxitch
\noindent
A 
logic
system 
$\alpha$
will be called 
{\bf Type-M} iff it contains
\eq{totdefxs} -- \eq{totdefxm}
as axioms,  
{\bf Type-A} iff it contains only
\eq{totdefxs} and \eq{totdefxa} as axioms,
and {\bf Type-S} iff it contains
only \eq{totdefxs} as an
 axiom. 
Also a  system is called 
{\bf Type-NS} iff it  contains
none of these axioms.
The 
significance of these constructs is explained by
items (a) and (b):
\bed
\item[  a.  ]
The existence of
Type-A systems that can  recognize 
their 
consistency under semantic tableaux deduction
was proven in \cite{ww5}.
Also, \cite{wwapal} demonstrated 
a large class of
Type-NS systems 
can 
recognize their
 Hilbert consistency.
(Many of
these systems 
do
prove
all 
Peano Arithmetic's
 $\Pi_1$ theorems
in a language
that represents addition and multiplication 
as 3-way predicates.)

\medskip

\item[   b.  ]
The above 
evasions of the Second Incompleteness
Theorem are known to be near-maximal in a mathematical sense.
This is because
the
combined work of Pudl\'{a}k, Solovay, Nelson and Wilkie-Paris
\cite{Ne86,Pu85,So94,WP87} implied  no
natural 
Type-S system can recognize  its  Hilbert consistency,
and   Willard
 \cite{ww2,ww7,ww9} 
strengthened 
some 
earlier
results by
Adamowicz-Zbierski 
\cite{Ad2,AZ1} 
to establish that  most
Type-M  systems cannot recognize their semantic
tableaux consistency.
\ennd

An unusual aspect of 
(a) and (b) is that there is a 
tight match between their 
positive and negative
results from a purely quantitative
perspective, but they
still
do not 
reach one's
ideal goals for 
deduction.
This is because there is something about
the 
deep yearnings
of research, 
roughly suggested by
Hilbert  in  $**$ and
 G\"{o}del in $*$, 
which has not 
been addressed.
This 
yearning and its
tiny remaining gap will be
our main focus  
in 
this article.

Other
efforts to 
evade the Second Incompleteness Theorem 
have used
the Kreisel-Takeuti    ``CFA''
system \cite{KT74}
or what 
could
be called
the interpretational framework of
Friedman,
Nelson, Pudl\'{a}k and Visser
\cite{Fr79b,Ne86,Pu85,Vi5}.
These methods
do not use
Kleene-like {\it ``I am consistent''} axioms,
similar  to
those in
 our work. 
Instead, CFA uses the 
unique
properties of Second Order generalizations of Gentzen's
Sequent Calculus (with modus ponens absent), and 
interpretational frameworks 
formalize how some systems 
recognize their
 Herbrand consistency 
on localized sets of integers
(which 
unbeknownst to 
themselves)
includes all 
natural numbers.

Such systems are not 
germane to our exploration
 of
Kleene-like {\it ``I am consistent''} 
axioms,
in the current article, 
but they do illustrate alternative
approaches that are germane to other
fascinating
open questions about 
Incompleteness paradigms.

As the reader exams the next several sections of this paper,
he should
also keep in mind that self-justifying systems that
rely upon semantic tableaux as their primary deduction method
is a different topic than formalisms using Herbrand-styled
deduction. This fact was brought to our attention
by
private communications from
L. A. Ko{\l}odziejczyk 
\cite{Ko5}, who noted that Herbrand styled systems can be exponentially
slower than semantic tableaux systems under 
well-defined
circumstances. 
 Ko{\l}odziejczyk's
insightful
 observation enabled us to prove in
\cite{ww9} 
that  Herbrand-styled self-justifying systems,
unlike their semantic tableaux counterparts, could house
a multiplication function symbol. 
Our results about mutliplication in 
\cite{ww9} 
are clearly of
mathematical interest,
but they are not germane to our philosophical
orientation in the current article because 
\cite{ww9}'s formalism is able to house a multiplication
function symbol only because of its exponential increase
in inefficency. 
Thus, the current article will focus instead on uses
of \cite{ww5,ww6}'s semantic tableaux 
styled
formalisms because we suspect that their added efficiencies
are likely to be of greater significance, from both a
philosophical and engineering perspective,
than our
slower and more inefficient
 alternative
mechanisms
 in 
\cite{ww9}.

\vspace*{- 0.8 em}

\section{Defining Notation}

\label{ss3-4}
\label{tst3}

\label{sect3}

The next two sections
will summarize
the properties of \cite{ww5}'s
IS$_D(A)$ axiom system. 
We
will
then outline how 
to  refine \cite{ww5}'s results
in the remainder of this article.

\cvu

A function
$ ~ \, F(a_1,a_2...a_j)~ \, $ 
will be called  
{\bf Non-Growth} 
iff  it satisfies the general inequality of
$ F(a_1,a_2,...a_j) 
\leq   Maximum(a_1,a_2,...a_j)$.
Six  examples of  
non-growth functions are
{\it Integer Subtraction} 
(where $~x-y~$ is defined to equal zero when
 $~x \leq y~),~~$
{\it Integer 
Division}
(where $~x \div y~$ is defined to equal
$~x~$ when $~y=0$, and
it equals $~\lfloor ~x/y ~\rfloor~$ otherwise),
$~~Maximum(x,y),~~$
$ Logarithm(x)$,$~$
$Root(x,y) \, =  \, \lceil  \, x^{1/y} \,  \rceil$ and
$Count(x,j)$  designating the number of ``1'' bits
among $ \, x$'s rightmost $ \, j \, $ bits.

\smallskip

The term
{\bf U-Grounding Function} 
referred in \cite{ww5} to
a set of eight operations, which included the
preceding
functions plus the {\it growth operations} of addition and
{\it Double$(x)=x+x$}. 
$\,$Our language $L^*$  was 
built
out of these 
function symbols, the usual
predicates of  ``$ \, = \, $'' and ``$ \, \leq \, $'' and 
the constant symbols
``0'' and ``1''.

\smallskip

In a context where  $~\, t \, ~$ is  any term in \cite{ww5}'s
language $L^*$, 
$\,$the quantifiers
stored in the wffs of
\newline
$ \forall ~ v \leq t~~ \Psi (v)~$ and $~ \exists ~ v \leq t~~ \Psi (v)$
were called {\it bounded 
quantifiers}.
Also, any formula in the 
\newline
U-Grounding language, all of whose
quantifiers are bounded, was called
a $\Delta_0^*$ formula.
The  $~\Pi_n^{* }~$ and  $~\Sigma_n^{* }~$ formulae 
were
then defined by the 
usual
rules that:
\bee
\small
\baselineskip = 1.1\normalbaselineskip
\xxitch
\item
Every  
$\Delta_0^*$ formula is considered to
be a
``$~\Pi_0^{* }~$''  and 
also to be a
 ``$~\Sigma_0^{* }~$'' 
formula.
\item
For $n \geq 1$,
 a formula 
is  called
 $ \,\Pi_n^{* } \,$
when it is encoded as 
$\forall v_1 ~ ...~ \forall v_k ~ \Phi$  with
$\Phi$ being  $\Sigma_{n-1}^{* }$
\item
Likewise,
 a formula 
is  called
 $~\Sigma_n^{* }~$
when it is encoded as 
$\forall v_1 ~ ...~ \forall v_k ~ \Phi,$  with
$\Phi$ being  $\Pi_{n-1}^{* }$.
\ene

\begin{example} \rm
 \label{exanew}
Although our language $L^*$ contains no multiplication
function symbol, it does provide a means to encode
multiplication as a 3-way relation 
$\Delta_0^*$ formula, Mult$(x,y,z)$, as 
is
illustrated
below:

\begin{equation}
\label{oneweq1}
\small
[~(x=0    \vee    y=0 ) \Rightarrow z=0~ ]~ ~\wedge ~~ 
[~(x \neq 0 \wedge y \neq 0~) ~ \Rightarrow ~
(~ \frac{z}{x}=y  ~\wedge \, ~  \frac{z-1}{x}<y~~)~]
\end{equation}
Moreover,
\eq{comm}
illustrates how the commutative principle for multiplication
can receive a $\Pi_1^*$ encoding
via the above  Mult$(x,y,z)$ predicate.
\begin{equation}
\label{comm}
\forall ~x~~\forall ~y~~\forall ~z~~~~\mbox{Mult}(x,y,z) ~
\Leftrightarrow~\mbox{Mult}(y,x,z) ~
\end{equation}
Similar methodologies
may also encode multiplication's
associative
and distributive axioms 
as $\Pi_1^*$ sentences.
In general, all the $\Pi_1$ theorems in a conventional arithmetic
language (that has multiplication symbols) can be
translated into equivalent $\Pi_1^*$ statements in our
language $L^*$. This implies that the set of
 $\Pi_1^*$ sentences 
is a quite rich class of statements and that no decision
procedure is capable of separating
all true from false
 $\Pi_1^*$ statements.
\end{example}

\bigskip

{\bf Further Terminology:} $~$Most of our other notation
will be similar to 
\cite{ww5}'s 
terminology.
In addition to considering a
definition of semantic tableaux deduction
that is similar to either
Fitting's
classic formalism
\cite{Fi90} or
a variant of it in 
\cite{ww5},
we 
considered in \cite{ww5}
also
a 
somewhat
 stronger technique, called
\txl{k} deduction.
It   is
defined below and
consists of a 
speeded-up version of a 
tableaux proof,
which 
permits an analog of modus ponens for 
the limited cases of performing Gentzen-style deductive
cuts on $\Pi_k^*$ and   $\Sigma_k^*$ formulae.

\smallskip

\begin{definition}
\label{newdef}
\rm
Let
 $~H~$ 
denotes a sequence of ordered pairs
$~(t_1,p_1),~(t_2,p_2),~...~(t_n,p_n),~$
where $~p_i~$ is a Semantic Tableaux proof of the theorem $~t_i.~$
Then  $H$ 
will be called a
{\bf Tab-k
Proof}
of a theorem $~T~$ 
from the axiom system $~\alpha~$
  iff $~T=t_n~$
and also:
\begin{enumerate}
\xxitch
\it
\item
Each of the ``intermediately derived theorems'' 
$~t_1,t_2, \, ... \, , t_{n-1}~$
must have a complexity no greater than that of
either a $\Pi_k^*$ or $\Sigma_k^*$ sentence.
\item
Each axiom in  $ p_i$'s
proof 
either
comes from $\alpha$ or is
 one of $ t_1,t_2, \, ... \, , t_{i-1} $. 
\end{enumerate}
$~~~~$ Thus, the preceding
2-part definition implies  a 
\txl{k} 
proof differs from 
Fitting's
definition 
of a
semantic 
tableaux 
proof
\cite{Fi90}
by allowing for
an analog
of a Gentzen cut-rule to be applied 
to intermediate results 
at  the
levels of  $\Pi_k^*$ and $\Sigma_k^*$ formulae.
\end{definition}

\medskip

\begin{remark} \rm
Let us say that 
an axiom system $\alpha$ 
has a {\bf Level-J Understanding}
of its own 
consistency 
under a 
particular
deduction method $D$ 
iff $\alpha$ can prove that there exists no proofs
using
its axioms and $D$'s deductive
methodology of both a
$\Pi_J^*$ theorem and its negation.
In this notation, items A and B summarize
\cite{sp0,ww2,wwlogos,ww5,ww7}'s results:
\bed
\xxitch
\item[   A.   ] 
 For 
every
axiom system $A$ using $L^*\,$'s
 U-Grounding language, 
\cite{ww5} 
showed its
IS$_D(A)$  formalism 
could prove
all $A$'s $\Pi_1^*$ theorems and simultaneously
verify its own Level-1
consistency when $D$ corresponds to
\txl{1} deduction.

\smallskip

\item[   B.   ] 
Two negative results that tightly complemented
item A's
positive result
were exhibited
in 
\cite{sp0,ww2,wwlogos,ww7}. The first
was that \cite{sp0,ww2,ww7} showed
most 
systems
are 
unable to verify their 
Level-0 consistency under
semantic tableaux 
deduction 
 when they included 
statement
\eq{totdefxm}'s ``Type-M''
axiom  that multiplication
is a total function. Moreover, \cite{wwlogos}
offered an alternate
form
of the
Second Incompleteness
Theorem
 that showed statement
\eq{totdefxa}'s
{\it 
far weaker} 
Type-A 
systems
are
unable to 
verify
their Level-0 consistency under
\txl{2} deduction.
\ennd
\end{remark}

The contrast between these
positive and negative results
has
 led to our conjecture that
sophisticated
theorem provers
are likely
 to 
eventually
achieve
a fragmentary part of the ambitions stated by 
G\"{o}del and Hilbert
in 
$*$ and $**\,$.
This is because
the question of whether a
formalism can support an 
{\it idealized Utopian}
conception of
its own consistency is {\it 
different} from 
exploring the degrees to which 
theorem-provers
can possess 
a {\it fragmentary
knowledge} of 
their own
consistency. 
Thus, the 
Incompleteness Theorem 
has demonstrated 
an Utopian idealized form of self-justification
is unobtainable,
but our 
on-going
research has found some
restrictive
forms of self-knowledge are, indeed,
feasible.

\section{The IS$_D(A)$ Axiom System}

\label{ss4-5}
\label{tst4}

\label{sect4}

In a context where
$~D~$ denotes any deductive method and
 $~A~$ denotes any axiom
system using
 $L^*\,$'s 
U-Grounding language, IS$_D(A)$ 
was defined 
in \cite{ww5}
to be an axiomatic
formalism  capable of recognizing all of 
$A$'s $\Pi_1^*$ theorems and 
corroborating 
its own Level-1 consistency 
under
deductive method  $D$.
It was  defined in
\cite{ww5} 
to consist of the 
following four
groups of axioms:
 
\begin{description}
\cvl
\item[Group-Zero:]
Two of the Group-zero axioms will 
be $\Pi_1^*$ statements defining
the
named constant-symbols,
$\bar{c}_0$
and $\bar{c}_1$,
 that designate the integers of 0 and 1.
The third and fourth Group-zero axioms
will be $\Pi_1^*$ statements defining the
two standard
growth functions of addition and 
$~Double(x) \, = \, x+x.~$
The
net effect of these 
axioms will be to set up a machinery to
define
any 
integer
$~n \geq 2~$ 
using fewer than
$3 \cdot \lceil \, $Log$~n~ \rceil \,$ 
logic symbols.

\medskip

\item[Group-1:] 
This axiom group  will consist of a
finite set of $\Pi_1^{*} $ sentences, denoted as $~F~$, which
which can prove any $\Delta_0^*$ sentence that
holds true under the standard model of the natural numbers.
(Any finite set of 
$\Pi_1^{*} $ sentences $~F~$ 
with this property
may be used to define Group-1,
as    \cite{ww5}  noted.)

\medskip

\item[Group-2:]
Let $\ulxyz \Phi \urxyz$ denote
$\Phi$'s G\"{o}del Number, and
HilbPrf$_A(\ulxyz \Phi \urxyz,p)$ denote a
$\Delta_0^{*} $ formula indicating 
$~p~$ is a
Hilbert-styled proof of 
theorem $~\Phi~$ from
axiom system $A$.
For each $\Pi_1^{*}  $ sentence  $\Phi$, the 
Group-2 schema will contain an axiom of 
form \eq{group2}.
(Thus IS$_D(A)$ can trivially prove
 all $A$'s 
$\Pi_1^{*}  $ theorems.) 
\begin{equation}
\forall ~p~~~\{~ \mbox{HilbPrf$_A(\ulxyz \Phi \urxyz,p)$}
 ~~
\Rightarrow ~~ \Phi~~\}
\label{group2}
\end{equation}
\item[Group-3:]
The final part of the IS$_D(\aaa)$ 
will consist of a single
self-referencing
$\Pi_1^*$
sentence, indicating 
IS$_D(\aaa)$ meets 
\textsection  
 \ref{tst3}'s criteria of being
``Level-1 consistent'' 
under deductive method $D$.
It 
is, 
thus,  the following  declaration:
\medskip
\begin{quote} 
\xxitch
\# $~$There exists no two
proof of
 both a $\Pi_1^{*} $ sentence
and its negation when
$D$'s deductive method is applied to an axiom system that
consists of  
the {\it union}
of the 
Group-0, 1 and 2 axioms with {\it this sentence
looking at itself}.
\medskip
\end{quote}
One encoding of \#,
$\,$as a self-referencing
$\Pi_1^{*} $ 
sentence, was provided in 
\cite{ww5}.
Thus, 
sentence \eq{group3} 
will be a
$\Pi_1^{*}\, $
encoding for  \#,
in a context where  
$ \, \mbox{Prf} \, _{\mbox{IS}_D(A)}(a,b) \, $ is
\nop
a 
 $\Delta_0^{*} $ formula 
indicating that $ \, b \, $ is a proof of a theorem $\, a\,$
via deduction method $D$
from  $\mbox{IS}_D(A)$'s axiom system,
and where 
Pair$(x,y)$ is a $\Delta_0^{*} $ formula
indicating  $ \, x \, $ is the G\"{o}del number of a
 $\Pi_1^{*} $ sentence and
 $ \, y \, $ represents $ \, x \,$'s negation.
\end{description}
\begin{equation}
\forall  ~x~\forall  ~y~\forall  ~p~\forall  ~q~~~~ \neg ~~
[~~ \mbox{Pair}(x,y)~ \wedge ~ 
~\mbox{Prf}~_{\mbox{IS}_D(\aaa)}(x,p)~
\wedge ~ ~\mbox{Prf}~_{\mbox{IS}_D(\aaa)}(y,q)~ ]
\label{group3}
\end{equation}

\cvl

{\bf Notation.} An operation $~I(~\bullet~)~$ that maps
an initial axiom system $\,\aaa \,$ onto an alternate
system  $\,I(\aaa)\, $ will be called {\bf Consistency Preserving}
iff  $\,I(\aaa)\, $ is consistent whenever all of $\aaa$'s axioms hold
true under the standard model of the natural numbers. In this
context, 
\cite{ww5} demonstrated:

\begin{theorem}
\label{ttt1}
\label{thold}
Suppose 
the symbol $D$ denotes either semantic
tableaux deduction or its \txl{1} generalization
(given in Definition \ref{newdef}).
Then the  IS$_D(~\bullet~)~$ mapping operation is consistency preserving
(e.g. 
IS$_D(\aaa) $ 
will be consistent whenever all of $\aaa$'s axioms hold
true under the standard model of the natural numbers).
\end{theorem}

We emphasize that the most difficult part of \cite{ww5}'s
result was 
neither the definition of its  
IS$_D(\aaa) $'s axiom system nor the
$\Pi_1^*$ fixed-point
 encoding of \eq{group3}'s Group-3 axiom.
Instead, it was the 
confirming
of \thx{thold} 
``Consistency Preservation''
property.

\smallskip

The 
confirming of
this
property
was 
subtle
because its invariant breaks down when 
$~D~$ is a deduction method only slightly stronger than
either semantic tableaux or
Definition \ref{newdef}'s
  \txl{1} construct. 
Thus, the combined work of Pudl\'{a}k and Solovay \cite{Pu85,So94}
implied that  \thx{thold}'s analog fails when $D$ represents
Hilbert deduction, and \cite{wwlogos} showed its generalization
would even fail when  $D$ represents  \txl{2} deduction.

These difficulties occur because a G\"{o}del-like
diagonalization effect will
automatically render
IS$_D(\aaa) $ inconsistent
when $D$ is  too strong.
Throughout this article, we will need to
assure 
$D$ does not
reach a strength that produces
such issues.

\section{A Finitized Generalization of  \thx{thold}'s Methodology}

\label{ss5-6}
\label{tst5}

\label{sect5}

One awkward aspect of
IS$_D(\aaa)$ 
is
that it
employs an infinite number of different
incarnations of
sentence \eq{group2}'s
Group-2 schema (since it contains one incarnation
of this sentence for each $\Pi_1^*$ sentence $\Phi$ in
$L^*\,$'s language). Such a Group-2 schema is 
cumbersome
because
it simulates $A$'s 
$\Pi_1^*$
knowledge almost via a brute-force 
enumeration.

Our Definition \ref{dd-is2} and Theorems
\ref{ttt2} and \ref{ttt3} will show how
to 
mostly
overcome this difficulty by  
compressing the infinite number
of
instances of sentence \eq{group2} in
IS$_D(\aaa)$'s Group-2 schema into
a purely finite structure.

\begin{definition}
\label{dd-is2}
\rm
Let $~\beta~$ denote any
finite set of
axioms that have 
 $\Pi_1^*$ encodings.
Then 
\I2 
will denote an axiom system,
similar to  IS$_D(\aaa)$, except 
its Group-2
scheme will employ $~\beta\,$'s set of axioms,
instead of using an infinite number of applications
of
statement \eq{group2}'s scheme.
(Thus,
the 
{\it ``I am consistent''} statement
in \I2's Group-3
axiom will be the same as before, except that
the  {\it ``I am''} 
fragment of its
self-referencing
statement
will reflect 
these
 changes in Group-2 in the obvious manner.)
\end{definition}

Our next theorem will
indicate
that an analog of \thx{ttt1}'s consistency preservation
property applies to \I2's
formalism:

\begin{theorem}
\label{ttt2}
Let 
 $D$ again denote either
semantic
tableaux deduction
or 
Definition \ref{newdef}'s
 \txl{1} construct,
and $\beta$ again denote a set of
$\Pi_1^*$ axioms.
Then
\I2 
will be consistent whenever all 
$\beta$'s axioms hold
true under the standard model.
(In other words, 
 \I2 
will satisfy an analog of \thx{ttt1}'s
consistency preservation property for IS$_D(\aaa) $.) 
\end{theorem}

A formal proof of \thx{ttt2} from first principles 
is
essentially
as lengthy as \cite{ww5}'s proof of 
\thx{ttt1}.
It is unnecessary to provide such a proof
here
 because 
both theorems 
are justified using 
roughly
analogous methods.
Appendix A 
offers a
brief
summary about how \cite{ww5}'s proof of  
\thx{ttt1} can be easily incrementally modified to
also  prove
\thx{ttt2}.

\smallskip

Our goal in the current article 
will be to show how
\thx{ttt2}'s
 \I2 formalism  
can assure that 
self-justifying axiom systems can prove
an astonishingly wide breadth of theorems,
{\it even when such logics}
contain only a
{\it strictly finite}
number of 
axiomatic statements.
We  need one further definition
to 
explore the depth of this meta-result:

\begin{definition}
\label{dkern}
\rm
Let
$~\ulcorner \Psi \urcorner ~$ denote
$\, \Psi\,$'s G\"{o}del number.
A $\Delta_0^*$ formula, 
Test$_i(t,x)$, 
will be called a {\bf Kernelized Formula}
iff Peano Arithmetic can prove every $\Pi_1^*$ sentence
$~\Psi~$ satisfies \eq{testker}'s
identity:
\beq
\label{testker}
\Psi ~~~ \Longleftrightarrow~~~ \forall ~x~~
\mbox{Test}_i  (~\ulxyz~\Psi~\urxyz~,~x~)
\enq
There are 
countably
infinitely many different 
 $\Delta_0^*$  predicates
Test$_1(t,x)$, Test$_2(t,x)$, Test$_3(t,x)$ ...
that satisfy this kernelized condition.
(One such kernel predicate is illustrated by
Example \ref{exa-1} below).
A formal enumerated list  of all 
such kernel
  predicates
will be called a  {\bf Kernel-List}.
\end{definition}

\begin{example} \rm
 \label{exa-1}
The set of
true $\Sigma_1^*$ sentences is
clearly
 recursively enumerable. 
This 
easily 
implies
there 
exists a $\Delta_0^*$ formula,
called say Probe$(g,x)$, 
such
that $~g~$ 
is
the G\"{o}del number of
a $\Sigma_1^*$ statement that holds true in the Standard
Model if and only if \eq{e-probe} is true: 
\beq
\label{e-probe}
\exists ~x~~~ \mbox{Probe}(g,x)~\wedge~ x \geq g
\enq
Now, let Pair$(t,g)$ denote a $\Delta_0^*$ formula
that specifies $~t~$ is the  G\"{o}del number of
a $\Pi_1^*$ statement and
 $~g~$ is
the  $\Sigma_1^*$ formula which is its negation.
Then our notation implies 
that
  $~t~$ 
is
a true 
 $\Pi_1^*$ statement 
if and only if \eq{e-2probe} holds true:
\beq
\label{e-2probe}
\forall ~x~~~ 
\neg~[~\exists ~g ~\leq~x~~~~~ \mbox{Pair}(t,g)~\wedge~\mbox{Probe}(g,x)~~]
\enq
Thus if
Test$_0(t,x)$
denotes the $\Delta_0^*$ formula of
``$~\neg~[~\exists ~g ~\leq~x~~ 
\mbox{Pair}(t,g)~\wedge~\mbox{Probe}(g,x)]~$'',
then 
it
is one example of what 
Definition \ref{dkern}
would
call a
``Kernelized Formula''.
\end{example}

\begin{definition}
\label{def3}
\rm
Let us recall that Definition \ref{dkern}
defined 
{\bf Kernel-List} to be an enumeration of
all the available kernelized formulae of 
Test$_1(t,x)$,
 Test$_2(t,x)$, Test$_3(t,x)$ ...
Then if Test$_i(t,x)$ is the $i-$th element in this
list of kernels and if $~\Psi~$ is an arbitrary $\Pi_1^*$ sentence,
then $\, \Psi \,$'s
{\bf i-th Kernel Image} will be  the 
following $\Pi_1^*$
sentence:
\beq
\label{imagker}
 \forall ~x~~
\mbox{Test}_{\, i \,} (~\ulxyz~\Psi~\urxyz~,~x~)
\enq
\end{definition}

\begin{example} \rm
 \label{exa-2}
The combination of Definitions 
\ref{dkern}
and \ref{def3} indicates that there is a
 subtle relationship
between a sentence $~\Psi~$ and its $i-$th Kernel Image.
This is because 
Definition \ref{dkern}
indicates that Peano Arithmetic does prove the invariant
\eq{testker},  indicating that 
 $~\Psi~$ 
is equivalent to its
 its $i-$th Kernel Image.
However, a weak axiom system 
can be plausibly uncertain about
whether this equivalence 
holds.

Thus if a weak axiom system proves statement 
\eq{imagker} (rather than $~\Psi~$ itself), 
then it 
may not be able to equate these two results.
This problem will apply to \thx{ttt3}'s
formalism.
However, \thx{ttt3} will be still of much interest
because \textsection  \ref{tst6} will 
illustrate a generalized methodology that overcomes
many of \thx{ttt3}'s limitations.
\end{example}

\begin{theorem}
\label{ttt3}
Let $~A~$ denote any
system, all of  whose
 axioms hold 
true 
in arithmetic's standard model,
and $~i~$ denote the index
of any of
Definition \ref{dkern}'s
kernelized formulae
 Test$_i(t,x)$.
Then it is possible to construct a
finite-sized 
collection of $\Pi_1^*$ sentences, called say
 $\beta_{A,i}$,
where 
\ik3
satisfies the following invariant:
\begin{quote}
\rm
If $~\Psi~$ is one of the 
$\Pi_1^*$ theorems of
 $~A~$
then  \ik3 can prove 
\eq{imagker}'s
statement 
 (e.g. it will prove the
``the $\, i-$th kernelized image'' 
of the sentence $\,\Psi\,$).
\end{quote}
\end{theorem}

\noindent
{\bf Proof Sketch:}
Our justification of 
\thx{ttt3} will 
use the following notation:
\bee
\item
Check$(t)$ will denote a $\Delta_0^*$ formula
that 
produces a Boolean value of True when
$t$ represents the G\"{o}del
number of a $\Pi_1^*$ sentence.
\smallskip
\item
 $~\mbox{HilbPrf}_A \,(   t   ,   q   )~$ 
will denote
 a  $\Delta_0^*$ formula that indicates
$~q~$ is a Hilbert-style proof of the theorem
$~t~$ from axiom system   $~A~$.

\smallskip

\item
For any kernelized
Test$_i(t,x)$ 
formula, the symbol GlobSim$_i$ 
will 
then 
denote \eq{globsim}'s $\Pi_1^*$ sentence.
(It will be called $A$'s $i-$th
{\bf ``Global Simulation Sentence''}.)
\ene
\beq
\label{globsim}
\forall ~t~~
\forall ~q~~
\forall ~x~~\{~~
[~~\mbox{HilbPrf}_A \,(   t   ,   q   )~~ \wedge ~~
\mbox{Check}(t)~~]~~~
\Longrightarrow ~~~ 
\mbox{Test}_i(t,x)~~~ \}
\enq

In this notation, 
the requirements of \thx{ttt3} will
be satisfied by any
version of the axiom system \I2, whose Group-2 schema $~\beta~$
is a finite sized
consistent set of $\Pi_1^*$ sentences
that has 
\eq{globsim}
as an axiom.
(This includes 
the minimal sized such system, 
that has only \eq{globsim} as an axiom.)
This is because if 
$\Psi$ is any 
$\Pi_1^*$ theorem of $A$, whose proof
is denoted as $~\bar{p}~$,  then both the
$\Delta_0^*$ predicates of
$\mbox{HilbPrf}_A \,( \ulxyz \Psi \urxyz , \bar{p}    )$ and
$\mbox{Check}(  \ulxyz \Psi \urxyz )$ 
are true. 
Moreover, \I2's 
Group-1 axiom subgroup was defined so that
it can automatically prove all
 $\Delta_0^*$ sentences that are true. 
Thus, 
 \I2 will
 prove these two statements and 
hence corroborate (via axiom
\eq{globsim}) the further statement:
\beq
\label{interm}
\forall ~x~~
\mbox{Test}_{\, i \,}(~  \ulxyz \Psi \urxyz ~,~x~ )
\enq 
Hence for each of the infinite number of $\Pi_1^*$
theorems that $~A~$ proves, the above defined
formalism will prove a matching statement
that corresponds to the 
 $\, i-$th kernelized image
of 
each
such  proven theorem.
 $~~\Box$

\section{Pragmatic L-Fold Generalizations of \thx{ttt3} } 

\label{ss6-7}
\label{tst6}

\thx{ttt3} 
is of
interest
because every axiom system $\,A\,$
will have
its formalism
\ik3 
prove the 
 $\, i-$th kernelized image of every
 $\Pi_1^*$  theorem that $A$ proves.
This fact is helpful
because 
\eq{testker}'s invariance
holds for all $\Pi_1^*$ sentences.
Moreover, our  
``U-Grounded''
$\Pi_1^*$ sentences
capture all 
Conventional Arithmetic's
{\it crucial}
$\Pi_1$ 
information
because they can 
view
multiplication as a 3-way 
 $\Delta_0^*$ 
predicate
Mult$(x,y,z)$
via 
\eq{neweq1}'s
encoding of this predicate.

\begin{equation}
\label{neweq1}
\small
[~(x=0    \vee    y=0 ) \Rightarrow z=0~ ]~ ~\wedge ~~ 
[~(x \neq 0 \wedge y \neq 0~) ~ \Rightarrow ~
(~ \frac{z}{x}=y  ~\wedge \, ~  \frac{z-1}{x}<y~~)~]
\end{equation}

One difficulty with the \I2
and \ik3 
systems
was mentioned by Example \ref{exa-2}.
It was that 
although Peano Arithmetic can corroborate
\eq{testker}'s invariance for every $\Pi_1^*$ sentence $\, \Psi \,$,
these latter two systems 
cannot  always also do so.

While there will
probably
never be a 
perfect method for
fully
 resolving this 
challenge,
there is
a  pragmatic engineering-style solution
that is often available.
This is essentially because
 our proof of \thx{ttt3}
employed
a formalism $~\beta~$ that 
used essentially
only one axiom sentence (e.g. 
\eq{globsim}'s $\Pi_1^*$ declaration ).

Since  the  \I2 formalism was intended 
for use
by any finite-sized system $\beta$,
it is clearly possible to
include any
finite number
of formally true
$\Pi_1^*$ sentences in $\beta$.
Thus for some fixed constant $L$, 
one can easily let
 $\beta$ include
$L$ copies of   \eq{globsim}'s axiom framework for a finite
number of different
Test$_1$, Test$_2$ ... Test$_L$ 
predicates, each of which satisfy Definition \ref{dkern}'s
criteria for being kernelized formulae. In this case, 
\I2 
will formally
 map each  initial 
$\Pi_1^*$ theorem 
$\Psi$
of some axiom system $~A~$ onto 
$L$ 
resulting
different
$\Pi_1^*$ 
theorems
of
the form
\eq{imagker}.

Our conjecture is
that
a goodly number of issues concerning
some logic-based
engineering applications,
called say $E$, will 
have convenient solutions via
self-justifying axiom systems
that follow the strategy 
outlined in
 the preceding paragraph.
Thus, we are suggesting that if 
 $~\beta~$ is a large-but-finite set of axioms, that
consists of 
$L$ copies of   \eq{globsim}'s axiom framework for 
different
Test$_1$... Test$_L$ 
predicates, then 
at least some
futuristic
engineering applications $~E~$
may possibly
 have their needs met by an
\I2 formalisms, when a
software
 engineer meticulously  chooses
\footnote{ \fgf 
In addition to having 
$L$ copies of   \eq{globsim}'s axiom framework, it will 
be
probably
 preferable (but not necessary)
to have some finite number of 
say
$M$ additional
general purpose $\Pi_1^*$ axioms  added to $\beta$, so as
to
enhance the power and  increase the efficiency 
of \I2.}
a proper
finite-sized
 $~\beta~$
that fits 
its
needs.

\medskip

\begin{remark} \rm
\label{rema-1}
The preceding discussion was not meant to overlook the fact
that the Second Incompleteness Theorem is a 
robust  result,
that will preclude self-justifying systems from recognizing
their own consistency under 
strong
definitions of 
consistency (as was
discussed by Sections \ref{ss2-3} and \ref{ss3-4}).
Our suggestion, however,  is
that computers
are becoming 
so 
powerful, in both speed and memory size as
the 21st century
is progressing, that   there 
 likely
will
emerge engineering-style applications
$~E~$
 that will benefit from  \I2's 
self-referencing
formalisms when a {\it large-but-finite-sized}
 $~\beta~$ is delicately chosen.
(Moreover, it is of interest to speculate
about whether such computers
can, at least partially, imitate a human's gut instinctive
{\it at least partial} confidence
\footnote{ \fgf While human beings presumably do not have a fully
robust conficence in their own consistency. the footnoted sentence refers
to the fact that they 
certainly
have a sufficient
{\it  partial} confidence of such to, 
at least,
reach the threshold needed for
motivating themselves to cogitate.}  
 in his own consistency.)
\end{remark}

\cvu

\begin{remark} \rm
\label{rema-2}
The next three sections will 
further
amplify on 
the preceding
themes.
They will make it unequivocally clear that the 
Second Incompleteness Theorem is too powerful a result
for any partial exception to 
overcome 
 its formalism,
in an entire and full sense.
 Yet at the same time, they
will suggest
that some aspects of the
statements $*$ and $**$, by G\"{o}del and Hilbert, did correctly
foresee that a 
logical system 
can, at least,
 bolster
 a  {\it partial limited-style
appreciation} of its own consistency.
\end{remark}

\section{Comparing the Properties of Type-M and Type-A
Formalisms} 

\label{ss7-8}
\label{tst7}

Let us recall that 
axioms
\eq{totdefxs}-\eq{totdefxm} indicated that Type-A 
systems
differ from Type-M formalisms by treating Multiplication
as a 3-way relation (rather than as a total function).
For the sake of accurately characterizing what our systems
can and cannot do, we have described our results as
being fringe-like exceptions to the Second Incompleteness
Theorem from the perspective of an Utopian view of Mathematics,
while  perhaps 
being
more significant results from an
engineering-style perspective of knowledge. Our goal in
this section will be to 
amplify
upon
this 
perspective by taking a closer look at Type-A and Type-M
formalisms.

\medskip

During our discussion, 
$   x_0,   x_1,   x_2,     ...    $ 
and  $   y_0,   y_1,   y_2,     ...  $
will denote the sequences defined by
the 
identifies of:

\vspace*{- 0.9 em}
{ \small
\baselineskip = 0.5\normalbaselineskip
\britch
\beq
\label{zs}
x_0~~ \, =~~ \, 2~~ \, =~~ \, y_0
\enq
\beq
x_{i+1}~~~=~~~x_{i}~+~x_{i}
\label{as}
\enq
\vspace*{- 1.0 em}
\beq
y_{i+1}~~~=~~~y_{i}~*~y_{i}
\label{bs}
\enq}

\vspace*{- 1.2 em}

\noindent
These sequences  represent the growth rates
that are produced, 
respectively,
by the axioms
\eq{totdefxa} and \eq{totdefxm}, 
that specify addition and multiplication are total functions.

Equations \eq{zs}-\eq{bs}
imply $ \, y_n \, = \, 2^{2^n} \, $ and
 $ \, x_n \, = \, 2^{ n+1} \,  \, $. 
Thus, the  $   y_0,   y_1,   y_2,     ...  $ sequence will
grow at a much faster rate than the
$   x_0,   x_1,   x_2,     ...    $ 
sequence. (This is because
 $ \, y_n \,$'s binary encoding will have an Log$(y_n)\, = \,2^n~$ 
length
while  $ \, x_n \,$'s binary encoding will have
a much smaller length 
of size Log$(x_n) \, = \, n+1.~)$

Our prior papers noted that the 
difference between these
growth rates 
was the 
reason that \cite{sp0,ww2,ww7}
showed 
all natural 
Type-M
systems, recognizing 
integer-multiplication as a total function, were unable
to recognize their 
tableaux-styled 
consistency ---
while \cite{ww93,ww1,ww5} showed 
some Type-A 
systems could 
simultaneously prove all Peano Arithmetic $\Pi_1^*$ theorems and
corroborate their own 
tableaux consistency.
Their gist
was that a
G\"{o}del-like diagonalization argument, which causes an
axiom system to become inconsistent as soon as it proves
a theorem affirming its own 
tableaux consistency,
stems from
the exponential growth in 
the  series  $   y_0,   y_1,   y_2,     ...$ .
(In other words, 
this growth
facilitates 
an intense
amount of self-referencing,
using 
the identity
Log$(y_n)\, \cong \,2^n~,~$
that will invoke the force of
G\"{o}del's seminal
diagonalization
machinery.)

These issues 
do raise the following question 
about proofs of
 the Second Incompleteness Theorems:
\begin{quote} 
\small
\baselineskip = 1.03\normalbaselineskip
$~\bullet~~$  How natural are exponentially growing sequences analogous
to  $   y_0,   y_1,   y_2,     ...  $, whose $n-$th member
requires 
more than
 $2^n$ bits to encode, $~$when such encodings
exceed the number of atoms in the universe if 
simply $~n\,> 100~$?
Does the  employment of
such a sequence, for corroborating the Second Incompleteness
Effect, suggest that this formalism is relying upon partially
artificial constructs? 
\end{quote}

\cvh

We will not attempt to derive a Yes-or-No answer to Question 
$~\bullet~$ 
because it
is one of those epistemological questions that can be
 debated
endlessly.
Our point is that  $~\bullet~$ 
probably does not require a definitive
positive or negative answer because both perspectives
are useful.
Thus,
the theoretical existence of a sequence  
integers 
of $   y_0,   y_1,   y_2,     ...  $, whose binary
encodings are doubling in length, is tempting
from the perspective of 
an Utopian view of mathematics, while 
awkward from an engineering styled 
perspective.
We therefore ask: {\it ``Why not be tolerant
of both perspectives? ''}

One virtue of 
this tolerance is 
it 
ushers in 
a greater understanding
of the statements $*$ and $**$ that G\"{o}del and
Hilbert made during the first third of the 20th century.
This 
is
because the 
Incompleteness Theorem
demonstrates 
no 
formalism can display
an understanding of its own consistency in an
idealized
 Utopian
sense. On the other hand,  
\textsection  \ref{tst6}
suggested
these 
remarks by G\"{o}del and Hilbert 
 might receive
more sympathetic interpretations, 
if one 
sought to explore
such questions from a less ambitious
almost engineering-style perspective.

Our
main  thesis is 
supported by a 
theorem
from \cite{ww6}.  It indicated that semantic tableaux
variations of self-justifying systems have no difficulty
in recognizing that an infinitized generalization of
a computer's
floating point multiplication (with rounding) is a total
function. The latter 
differs from integer-multiplication,
by not having its output become double the length of
its input when a number is multiplied by itself.
Thus, the 
intuitive
reason that 
\cite{ww6}'s
 multiplication-with-rounding operation
is compatible with self-justification is
because it
 avoids the 
inexorable
exponential
growth, associated with \ep{bs}'s sequence 
 $   y_0,   y_1,   y_2,     ... ~  $.

Also, the \thx{ttt4} below indicates our self-justifying system
can recognize a double-precision variant of integer multiplication
as a  total function.

\begin{theorem}
\label{ttt4}
Consider the version of the IS$_D(\aaa)$ and \ik3 axiom systems
where $A$ denotes Peano Arithmetic. These systems can formalize
two  total functions, called Left$(a,b)$
and Right$(a,b)$, 
where any ordered pair $(a,b)$
is mapped onto
the objects,
 Left$(a,b)$
and Right$(a,b)$, designating
 the  bit-sequences
that represent
the left and right halves of
 the multiplicative
product of  $a$ and $b$.
\end{theorem}

\thx{ttt4}
follows from a straightforward generalization of
\cite{ww6}'s
analysis of
floating point multiplication.
Both
it 
and 
the \cite{ww6}'s result
suggest that
many arithmetic operations, appearing
in an engineering environment, are
nicely compatible with self justification.
This is because most of the engineering and
applied-mathematics uses of multiplication
can replace a purist variant of integer multiplication
with either
floating point 
multiplication
or \thx{ttt4}'s notion of double-precision
multiplication.

\section{A Different Type of Evidence Supporting 
Our
Thesis}

\label{ss8-9}
\label{tst8}

Let us recall that 
 Pudl\'{a}k and Solovay 
\cite{Pu85,So94} 
observed
that 
essentially all
Type-S
systems, 
containing merely
statement  \eq{totdefxs}'s
axiom that successor is a total function,
are unable to verify their own consistency under 
Hilbert deduction. 
(See also related work by 
Buss-Ignjatovic \cite{BI95},
 \v{S}vejdar \cite{Sv7} 
and 
the
Appendix A of \cite{ww1})

It turns out that
\cite{wwlogos} generalized 
these
 results to
show that Type-A 
systems are unable to verify their
own consistency, under Definition \ref{newdef}'s 
\txl{2} deduction methodology.
At the same time,  
Theorems \ref{ttt1} and \ref{ttt3} 
demonstrate that 
the IS$_D$
and IS$^{\#}_D$
axiomatic frameworks can verify
their own consistency under
\txl{1} deduction. Our goal in this section will be to
illustrate how the contrast between these positive and negative
results 
is
analogous to the differing growth rates
of the 
sequences
$   x_0,   x_1,   x_2,     ...    $ 
and  $   y_0,   y_1,   y_2,     ...  $ 
from Equations
\eq{zs}--\eq{bs}.

During our discussion 
$~G_i(v)~$ will denote 
the scalar-multiplication
operation  that maps
an integer $~v~$ onto 
$~ 2^{2^i}\cdot v~$. 
Also, $~\Upsilon_i~$ will  denote 
the statement, in the U-Grounding language, that 
declares that 
 $~G_i~$ is a total function.
Our paper \cite{wwlogos}
had noted that  $~\Upsilon_i~$ can receive
a $\Pi_2^*$ encoding. It is also obvious that $~G_i~$
satisfies the  identity:
\beq
\footnotesize
\britch
\label{e-Gi}
G_{i+1}(v) ~~~ = ~~~ G_i(~ \, G_i(v)~ \, )
\enq
It was 
noted in \cite{wwlogos} that 
this identity
implies  one
can construct 
an axiom system $  \beta  $, comprised of
solely $\Pi_1^*$ sentences,
where
a semantic tableaux proof 
can establish
$  \Upsilon_{i+1}$  
from
$  \beta+\Upsilon_i$  
in a constant number of steps. 
This implies, in turn, that a \txl{2} proof from
$  \beta  $ will require no more that O$(n)$ steps
to prove $  \Upsilon_{n}$ (when it uses the obvious
n-step process to
confirm in chronological order 
$~\Upsilon_1 \, , \,  \Upsilon_2 \, , \,  ... \Upsilon_n ~.~~)$

\smallskip

These observations are significant because 
$G_n(1)=2^{2^n}$.
Thus,
\cite{wwlogos} showed a \txl{2} proof
from $\beta$ can verify
in  O$(n)$ steps   
that this integer exists.

\smallskip

This example is  helpful because it illustrates
that the difference between the growth rates, 
under
proofs
using Definition \ref{newdef}'s 
 \txl{1} and \txl{2} deductive methodologies, is  analogous
to the 
differing
growth speeds
of
the
sequences $   x_0,   x_1,   x_2,     ...    $ 
and  $   y_0,   y_1,   y_2,     ...  $ 
from Equations
\eq{zs}--\eq{bs}.
Thus once again, the faster growth-rate has the side-effect
of triggering off the force of the Second Incompleteness Theorem
(e.g. see \cite{wwlogos}).

This analogy suggests
that the 
Second 
Incompleteness
Theorem has different implications from the perspectives
of 
Utopian and engineering
theories about
 the intended
applications of mathematics. Thus, a Utopian
may  possibly be 
 comfortable
with 
a
perspective, that contemplates sequences
 $   y_0,   y_1,   y_2,     ...  $ 
with
elements growing in length
at an exponential speed, but many engineers may be
suspicious of such
growths.

A hard-core engineer,
in contrast, might
 surmise that the inability of self-justifying
formalisms to be compatible with
Definition \ref{newdef}'s 
 \txl{2} deduction is
not 
as disturbing
 as it might 
initially 
appear to be.
This is
because \txl{2}
differs from 
 \txl{1} deduction
by producing 
growth rates
that are so sharp
that their material realization has no analog
in the everyday mechanical reality that is the
focus of an engineer's 
interest.

Our personal preference is for
a perspective lying
half-way
between 
that of an Utopian mathematician and
a hard-nosed engineer. 
Its
dualistic
approach
suggests
some form of half-way agreement
with the the goals of
Hilbert's
consistency program in  $**$.
(This is
 the maximal type of agreement with
Hilbert's goals that is, 
obviously,
feasible 
because
no more than a curtailed form of 
the
objectives
of his consistency program
is plainly
realistic.)

\section{Further
Related Improvements}

\label{ss9-10}
\label{tst9}

An added point is that there are many 
types of
self-justifying systems  available, with some
being better suited for  engineering environments
than others.

For instance, our initial 1993 paper \cite{ww93}
employed a Group-3 {\it ``I am consistent''} axiom
that was much weaker than 
the current  specimen.
The distinction was that
\cite{ww93}'s self-consistency declaration 
excluded 
merely
the existence of a semantic tableaux proof
of $0=1$ from itself, while 
the
sentence \eq{group3} is
more elaborate because
it excludes the existence of simultaneous proofs
of a $\Pi_1^*$ theorem and its negation.

This distinction is 
significant because 
it implies that each time our
self-justifying systems $~S~$  prove a $\Pi_1^*$
or  $\Sigma_1^*$ theorem $~\Phi~$, they will know
it is  fruitless
 to search for a
proof 
of the negation
of $\Phi$'s statement
\footnote{ \fgf From at least an engineer's
hard-nosed
perspective, this fact is
very helpful
because it will allow an automated
theorem prover to imitate, 
for at least
$\Pi_1^*$
and  $\Sigma_1^*$ sentences, 
the common-sense
human presumption that 
the 
discovery of a 
proof a theorem $\Phi$ 
leads to the time-saving conclusion that it is unnecessary
to search for a proof of $\neg ~\Phi$.}.

\def\britch{\small}

Ideally, one would  like to
additionally  develop self-justifying
systems $~S~$ that  could corroborate the validity
of \eq{brxefl}'s reflection principle for all sentences 
$\Phi$, as well.
\beq
\label{brxefl}
\forall p ~~[~ Prf_S^D(\ulxyz \Phi \urxyz,p)
  ~~ \Rightarrow  ~~ \Phi~~]
\enq
L\"{o}b's Theorem 
establishes,
however,
 that all
axiom systems $S$,
containing the strength of  Peano Arithmetic, are able to prove
\eq{brxefl}'s invariant 
{\it only in the degenerate case} where they prove $\Phi$
itself. Also, the Theorem 7.2 from \cite{ww1}
showed 
essentially all
axiom systems,
{\it weaker} than Peano Arithmetic, are unable to prove \eq{brxefl}
for all $\Pi_1^*$ sentences $\Phi$
simultaneously. 
\thx{ttt5}'s reflection principle
is thus 
pleasing:

\begin{theorem}
\label{ttt5}
For any
input axiom system $A$,
it is possible to extend the self-justifying
IS$_D(\aaa)$ and \ik3
systems,
from  Theorems \ref{ttt1} and \ref{ttt3}, 
so
that a new broader system 
$S$
has all the properties of the prior systems and can
additionally:
\end{theorem}

\vspace*{- 1.2 em}
\bee
\cvz
\it
\item
Verify that \txl{1} deduction supports the
following analog of 
\eq{brxefl}'s 
self-reflection principle
under $S$
for any
$\Delta_0^*$ and $\Sigma_1^*$ 
sentences $\Phi~~$:
\beq
\small
\footnotesize
\britch 
\label{nrxefl}
\forall p ~~[~ Prf_S^{\rm Tab-1}
(\ulxyz \Phi \urxyz,p)
  ~~ \Rightarrow  ~~ \Phi~~]
\enq
\item
Verify 
\eq{rdilute}'s more general
{\bf ``root-diluted''} reflection principle
for $~S~$ 
whenever
$\theta$ is $\Sigma \, _{1}^*$ 
and
 $\Phi$ is a $\Pi_2^*$ sentence of the
form ``$~\forall u_1  ... \forall u_n~~    
  \theta(u_1... u_n  )~$''. 
\beq
\label{rdilute}
\small
\britch 
\footnotesize
\forall p ~~[~ Prf_S^{\rm Tab-1}
(\ulxyz \Phi \urxyz,p)
  ~~~ \Rightarrow  ~~~ \forall x~~~
 \forall u_1< \sqrt{x}~~~    ... ~~~\forall u_n< \sqrt{x}~~~~    
  \theta(u_1... u_n  )~~]
\enq
\ene

\cvh

\thx{ttt5} 
is stronger than
Theorems \ref{ttt1} and
 \ref{ttt3} because of
its  reflection principles.
Its proof is
a generalization of related reflection princibles used
in \cite{ww1}. It is
 summarized 
in 
Appendix B.
\thx{ttt5} 
reinforces
our theme about how 
exceptions 
to
the Second Incompleteness Theorem 
may 
appear to
be 
minor 
from the perspective of
an Utopian 
view of mathematics,  
while being
significant
from an  engineering standpoint.
This is because:
\bed
\parskip 0pt
\item[A. ]
The ability of  \thx{ttt5}'s
system $S$ to 
support
\eq{nrxefl}'s 
self-reflection principle
under
\txl{1} 
proofs for
arbitrary
 $\Delta_0^*$ and $\Sigma_1^*$ sentences, 
as well as
to support
\eq{rdilute}'s
 root-diluted reflection principle 
for  $\Pi_2^*$ sentences,
is 
clearly
non-trivial.

\smallskip

\item[B. ] 
The incompleteness result
in  \cite{ww1}'s
Theorem 7.2, shows,
however, that the above reflection
principle cannot be extended to all
  $\Pi_1^*$ sentences, 
{\it in an undiluted sense.}
\ennd
The tight fit 
between
 A and B
is
analogous to
several other  
slender borderlines
that separated
generalizations and boundary-case exceptions
for the Second Incompleteness Theorem, that we had mentioned earlier.
The Second Incompleteness
Theorem 
should,
thus, 
 be 
seen
 as robust,
from an 
idealized
Utopian perspective on mathematics,
while 
permitting
caveats
from 
possibly
an
engineering 
styled 
perspective.
This
 dualistic perspective
allows one to
nicely
share at least some
{\it   partial
agreement} with the general spirit of the
 the
statements
$*$ and $**$ 
by  G\"{o}del and
Hilbert,
while simultaneously
 appreciating
the 
 stunning
achievement
of
the Second Incompleteness Theorem.

\section{``Wir m\"{u}ssen wissen: wir werden wissen''}

\label{newwir}

In summary,
our research has been
largely
 motivated
by the 
approximate
conjecture that 
human beings
 acquire the needed
mental energy and will-power for
motivating their
time-consuming
cogitations 
{\it only via}
 the human mind
owning, {\it at least in a weak sense}, 
some type of
 quasi-automatic and possibly unconscious
appreciation of its own
consistency.

Two of Hilbert's famous 
often-quoted statements suggest
that he would
agree, at least partially,
 with the main spirit of our conjectures.
The first was his 1925 statement $**$ 
(whose relationship to G\"{o}del's
supportive remarks in $*$
was already
noted
 in
\textsection   \ref{sss1} ).
The 
second 
Hilbert
statement
was
the widely known
motto of his consistency program,
``Wir m\"{u}ssen wissen: wir werden wissen'',
whose
English translation
appears
below. 
Part of what makes $***$ 
intriguing
 is that
Hilbert arranged
for this 
motto, 
often cited as the justification of
his consistency program,
to be placed on his tombstone, even though
he co-authored with Bernays  an historic generalization of the
Second Incompleteness Theorem, using their
     ubiquitous
 Hilbert-Bernays
Derivability Conditions \cite{HB39,Mend}. 
\begin{quote}
$***~$ ``We must know: we will know''
\end{quote}
The presumable reason that Hilbert believed that,
at least, some diluted version
of his consistency program would
ultimately succeed
 is because it
adamantly
 defies
common sense to explain how a human mind can motivate itself
to cogitate, without {\it possessing some at least diluted} form of knowledge
and faith in its own consistency.

Thus, Hilbert clinged to the belief that some 
{\it at least diluted version} of his  consistency
program would ultimately succeed, as $***$ had
hinted
\footnote{ \fgf\label{xcare}
Some readers may initially suspect that the 
``Wir m\"{u}ssen wissen: wir werden wissen'' statement was inscribed
on Hilbert's tombstone to make a general statement about his
philosophy of life, rather than to address remaining issues
left open by 
 G\"{o}del's Second
Incompleteness Theorem. While this point may be partially
true, the significance of $***$ being inscribed on his tombstone
seems to transcend 
this viewpoint.
This is because 
Hilbert
was
famous for having
adamantly
 continued,
throughout his life,  seeking
 a ``revised'' version of
his consistency program, even when he
co-authored 
\cite{HB39}'s 
generalization of  the
Second Incompleteness Theorem.
Thus, it appears to be no coincidence
that Hilbert, who had a knack for 
formulating
many famous often-quoted phrases,
would have the motto of his never-terminated
consistency program inscribed
as the
final
 memorial
 on his tombstone.}.
Moreover, a word-by-word reading of $***$
finds it 
openly
declaring some
mysterious
object XYZ 
``must'' exist,
whose exact nature
(and even defining name)
 it is
deliberately
vague
about (despite the fact that 
Hilbert tells us that ``we must know'' about it).

The prior nine sections of
the current article,
as well as our earlier research 
in \cite{ww93}-\cite{ww11}, were
largely devoted
to defining this
elusive object XYZ (that underlines a human's instinctive
faith in his own processes
but which is
decisively
 awkward to describe).
Part of the reason that
Hilbert's desired
 XYZ is so 
tenderly
elusive is
that \textsection\ref{ss7-8} 
has documented that systems, recognizing  {\it even}
their
 own 
{\it tiny}
 miniaturized
tableaux consistency, cannot corroborate the assumption
that multiplication is a total function, and 
\textsection\ref{ss8-9} has formalized how 
Type-A axiom
 systems are also
unable to recognize their own consistency  
under a 
\txl{2} generalization of semantic tableaux deduction.

These facts certainly 
illustrate
 that a formalism cannot recognize its
own consistency in an idealized Utopian sense,
and one must thus approach this subject matter
with 
daring and
ginger
amounts of caution.

Within this context, the gist of our Sections \ref{ss7-8} and \ref{ss8-9}
was that Type-A systems that recognize their own consistency
under
\txl{1} deduction are not,
actually, quite as weak as they
might initially appear.
This is because they avoid the super-exponential growth
of 
\ep{bs}'s 
  $   y_0,   y_1,   y_2,     ... ~ $ series,
whose growth is inscrutable from a hard-nosed engineering style
perspective. Thus, if one confines one's attention
onto
\ep{as}'s
more 
modest 
exponential
growth
 for its
sequence 
  $   x_0,   x_1,   x_2,     ... ~ $, then
self-justifying
 systems can be formalized
that 
are compatible with 
this slower growth 
rate
and
meet at least a 
diluted version of $***$'s objectives.
(These weak-styled systems, obviously, support only a
{\it 
tenderly
modest} knowledge of their own 
self-consistency,
relative to an idealistic Utopian's
aspirations for what mathematics should
ideally
 embody.
Their outlook
about self-consistency 
are, however, 
{\it not entirely} 
inconsequential
from an engineer's 
applications-oriented 
perspective.)

The thesis, advanced in this article, is 
subtle 
because there are
many different kinds of {\it ``I am consistent''} statements that
can be used by the Group-3 components of our self-justifying axiom
systems. Some of these formalisms have the advantages of being
preferable from a pedagogic perspective.
This is because a formalism that merely declares the
non-existence of a semantic tableaux proof
of $0=1$ from itself is 
relatively easy  to analyze
(such as the initial
1993
 version of a self-justifying axiom
system that was proposed in \cite{ww93}).
Other variants of the
self-justification
 construct,
which assert the non-existence of 
simultaneous \txl{1} styled
 proofs of a $\Pi_1^*$ sentence and its negation, 
are 
substantially
more
complex to analyze but have the advantage of 
providing
greater
levels of self-knowledge.
This second topic,
obviously, 
comes
closer
to meeting  Hilbert's 
hopes,
subsequent to 1931,
because 
its 
self-justifying formalisms 
contain
 a
larger degree of
inner strength and self-understanding.

Finally,
it is useful to conclude our 
perspective about 
a limited
partial revival
of 
Hilbert's 
consistency program
 by reviewing the contrast between
Theorems \ref{ttt1} and \ref{ttt3}. These theorems formalize a
fundamental trade-off, where each result has 
sharply
different advantages.
Thus, \thx{ttt1} indicates that every consistent axiom system $~A~$
can be mapped onto a recursively enumerable 
self-justifying
axiom system
IS$_D(A)$ 
that can formalize all $A$'s $\Pi_1^*$ knowledge 
while recognizing its own Level-1 consistency. 
In contrast,
\thx{ttt3} shows that the infinite number of different proper axioms
in IS$_D(A)$ can be reduced to a finite size if one 
is willing to settle for
 a system 
\ik3 that
replaces a full knowledge of all A's
 $\Pi_1^*$ theorems with 
 an $i-$th kernelized knowledge of this information.
Since part of Hilbert's goal was to find a means whereby systems
{\it of purely finite
size} can recognize their own consistency,
we suspect that \thx{ttt3}, along with Section \ref{ss6-7}'s enhancement
of its results, 
demonstrate
that some
fragmentary 
(but
not full)
part of Hilbert's
objectives in $**$ can be achieved.

\section{Miniaturized Finitism}
\label{mini}

The phrase {\it ``Miniaturized Finitism''}
should, perhaps,
be employed to describe the type of 
agenda for
mathematics
that this article 
is advocating.
The 
cautious-sounding
adjective
of {\it ``miniaturized''} was
attached to this name because one should 
readily 
admit
 that
there is an
indisputably
 miniaturized and 
quite
 humble
aspect to
any
school of
mathematics that drops the assumption that multiplication
is a total function and  which also replaces an
unqualified modus ponens
rule
 with Definition \ref{newdef}'s more tender
\txl{1} variant of deduction.

The second 
word 
 {\it ``Finitism''},
on the other hand,  was attached to 
this name 
because 
it
 partially vindicates the 
``finistic''
presumption
of
G\"{o}del and Hilbert, specified in
statements $*$ and $**$, that it should be 
feasible to
formalize how Thinking Beings 
do conceptualize  at least a partial {\it intuitive} appreciation
of their own consistency.
Thus, Miniaturized Finitism looks towards
the   
remaining 
approximate 
10 \% of the
issues 
raised by Hilbert's
Second Open Question, in a context where the Second Incompleteness
Theorem has resolved 
the more main-stream
$ 90 \, \%$ 
of the issues, raised by Hilbert's
 year-1900
open question, in a decisively negative 
manner \footnote{ \fgf 
Researchers who happen to
examine Willard's research since 1990 will
curiously notice
that the philosophy of the footnoted sentence has been
central to 
this
research since 1990. Thus, Fredman-Willard
\cite{FW90,FW93}
published in that year a boundary-case exception to
the widely accepted textbook lower bound for 
sorting and searching that would later become the chronologically
first among six items listed in the
{\it Mathematics and Computer Science} section of the
{\it 
National Science Foundation's 1991 Annual Report.}
And so
beginning with \cite{ww93}'s 
initial
1993 theorem about self verification,
Willard's on-going
 research
into Logic
was based on the
presumption that other types of widely-accepted 
enticing
textbook
results would be found, upon 
careful
inspection, 
ultimately
to 
admit partial, 
 {\it although certainly weak,}
forms of
{\it very tightly} defined boundary-case exceptions.}.

Some may
critically
 view 
our
particular
 form of self-justification
 as 
an essentially 
dwarfed-sized theory,
which 
uses an argument resting on a
tiny
 ``miniaturized'' notion of growth
to explain how Thinking Beings are
able to muster an at least partial appreciation of their
own consistency.
In such a context,
it is 
pleasing to reply that G\"{o}del was 
famous
for 
telling
Einstein that he 
did dearly
{\it love} dwarfs: For instance,
Dawson \cite{Da97} and
Goldstein \cite{Go5} recall
how
G\"{o}del 
``was espeically fond of Disney films''
and 
attended ``at least three''
repeated
screenings of
 ``Snow White and the Seven Dwarfs''.
Likewise,
Yourgrau \cite{Yo5}
recites 
 G\"{o}del
telling
Einstein 
 that
the tale of the  seven dwarfs 
was, 
actually,
 his favorite movie
 because
it
``presents the world as it should be and as if it had
meaning''.

Leaving
aside the 
poetry of the
preceding
metaphor about how G\"{o}del was enamored by
miniaturized-sized 
objects such as dwarfs,
our main
points are that both
G\"{o}del and Hilbert specified in
statements $*\,$ and $**$ that they 
felt some type of self-justification would be feasible.
Thus,
the  name {\it ``Miniaturized Finitism''},
 has
its adjective 
 {\it ``miniaturized''}
nicely capturing 
 the fact that
only a discernible
{\it fractional} 
(but 
indeed
quite non-trivial) part 
of the
initial goals of
Hilbert's consistency program
are
achievable.

\section{Summary Statement}

\label{conc}

As a final note, it is useful to recall that the opening chapter
of this article indicated the current
report
would have
different goals than our earlier papers.

The distinction was that our earlier work
consisted of papers employing mostly a mathematical gendre,
focusing on generalizations and boundary-case exceptions for the
Second Incompleteness Theorem. The emphasis in the current 
report
was different because it sought to focus on interpreting the meaning
of the statements $*$, $**$ and $***$ of G\"{o}del and Hilbert,
rather than centering around
the mathematics of our new formalisms.
(Our current paper presented new
theorems only when they were necessary to formalize
its epistemological overview.)

The reason for the current article's different 
emphasis
 is that
we anticipate
that this second topic
should
grow 
in
significance, as the 21st century progresses. This is 
mostly
because
 digital computers have grown
in importance during the
last seventy years and have essentially doubled in memory size every 2 years
at the unabated exponential rate that Moore's Law 
Moore's Law has predicted.

This fact about the burgeoning growth in memory size for digital
computers in turn raises the questions about their other capacities.
In particular, it raises the question of whether such computers
will have the ability to in some  sense
to recognize their own consistency, which ability humans seem
to possess in some partial, though not full sense.

Regarding the last point,
our speculations are
 not 
intended
to suggest that humans
(or computers) will be
able to
corroborate
their own consistency in a robust sense.
It is obvious that  
the
Second Incompleteness Theorem has 
established
that such a prospect is hopeless. The point is, however,
that humans have been able to intuitively appreciate their own
consistency in some type of quasi-reasonable sense --- in order
to motivate themselves to gain the mental energy and will power
to cogitate.

The
advance of
the digital computer during the 21st century 
will, thus,
 inevitably raise
the question about how computers can  maintain,
likewise, some type
of timid-but-well-defined appreciation of their own consistency.
We suspect this will cause 
the statements $*$, $**$ and $***$ of G\"{o}del and Hilbert
to be
reviewed often,
again, in the future
(and 
ultimately 
to gain some form of partial-although-not-full
 reinforcement).

This point is 
delicately
subtle because
there is no question that the Second Incompleteness
Theorem will always be remembered as the greatest
of all 
the
achievements in the
modern
 advancement of  Logic.
Yet, this fact should not preclude
the  Second Incompleteness Theorem 
as being seen as 
analogous
to 
many other
mathematical discoveries, in that 
most great mathematical results
typically
permit for the existence of
some types of 
well-defined 
forms of
 boundary-case exceptions.

Thus along the lines that
G\"{o}del and Hilbert had
approximately
 predicted, in 
especially 
their first and third
 statements
of $*$ and $***$,
it is not surprising that there exists certain
types of
unconventionally framed
axiom systems that possess
a form of
 diluted
{\it but
not fully
immaterial}
knowledge about their own self-consistency.

\bigskip

  {\bf Acknowledgments:} $~$I thank 
  Bradley Armour-Garb  and Seth Chaiken  for several
  useful suggestions about how to improve the presentation.

\newpage

\section*{Appendix A  summarizing  \thx{ttt2}'s Proof}

The same methods that 
Section 5 of
\cite{ww5} had used to prove
\thx{ttt1} 
can be easily generalized to
corroborate \thx{ttt2}'s similar consistency preservation
result. This is because   
\thx{ttt2}'s hypothesis implies
that all the Group 0, 1 and 2 axioms of
\I2 
are analogous to their 
counterparts
 under IS$_D(A)$,
in that both represent
$\Pi_1^*$ sentences  holding true under the
standard model of the Natural Numbers.
In this context, a similar reductionist argument, as had appeared
in Section 5 of
\cite{ww5}, will imply 
\I2's Group-3 axiom 
must also hold true.
Thus,
a direct
 analog of
\cite{ww5}'s syllogism 
will show that it is impossible
for an ordered pair $(p,q)$ to
contradict
 \I2's
 Group-3 axiom by 
simultaneously
having 
1) $p$ prove a $\Pi^*$ theorem, 2) $q$ prove its
negation and 
3)
Max$(p,q)$ represent the minimal
value among all ordered pairs with this property.

We will not delve into further details here because
\thx{ttt2}'s proof is
similar 
to
\cite{ww5}'s analogous
theorem  in that they both, essentially, 
rest on the fact that the Group 
 0, 1 and 2 axioms of their formalism
represent true $\Pi_1^*$ sentences. 
The reason for our interest in 
\thx{ttt2} in the current article is that it
is a useful
intermediate step for establishing
\thx{ttt3}. 
The latter,
in turn,   will  
reinforce the epistemological
viewpoint of
Sections \ref{ss6-7} --
\ref{conc}
(especially as it
 pertains to 
the question about whether a 
limited subset
 of
Hilbert's ``finistic goals''
can 
be revived in a
 fragmentized sense).

\section*{Appendix B summarizing  \thx{ttt5}'s Proof :}

Let $\Phi$ denote a prenex sentence of the form
$\forall u_1~\exists w_1~\forall u_2~\exists w_2...
\forall u_n~\exists w_n~~\theta(u_1,w_1...u_n,w_n)$
where $\theta(u_1,w_1...u_n,w_n)~$ is a $\Delta_0^*$
formula. Also, let us
define
$~\Phi^\oplus~$ 
to be the 
following
sentence:
\beq
\label{tangf}
\footnotesize
\britch
\forall x ~~~ \forall u_1<\sqrt{x} ~~~
\exists w_1~~~ \forall u_2~<\sqrt{x} ~~~
\exists w_2~~.....~~\forall u_n<\sqrt{x} ~~~
\exists w_n~~~
\theta(u_1,w_1...u_n,w_n)~~~
\enq
This notation was designed so that the 
operation that maps $\Phi$ onto 
$\Phi^\oplus$ will  change 
the ranges 
{\it of only}
 the unbounded
universal quantifiers in 
$\Phi$
(e.g. all
{\it bounded} universal quantifiers
appearing inside 
its
expression $\theta$
are left unchanged). 
Thus, $~\Phi ~ =~\Phi^\oplus~$ 
in the degenerate
case where  $~\Phi~$  is either a $\Delta_0^*$ or $\Sigma_1^*$
sentence.

An axiom system $S$ will be said
to own
 a
 ``tangible understanding'' of its
reflection
properties under deduction method D iff it can
corroborate
\eq{tangref} 
for every sentence $\Phi$.
\beq 
\footnotesize
\britch
\label{tangref}
\forall p ~~[~ Prf_S^D(\ulxyz  \Phi \urxyz ,p)
  ~~ \Rightarrow  ~~ \Phi^\oplus~~]
\enq 
The Theorem 6.1 of \cite{ww1} 
showed that
every 
consistent  axiom system $A$ can be mapped onto
a system $S$ that can prove all
$A$'s $\Pi_1^*$ theorems,  while verifying \eq{tangref}'s
reflection principle when $~D~$ denotes semantic tableaux deduction.

Moreover,
it is possible to
hybridize 
\cite{ww1}'s 
result with
\thx{ttt3}'s methodology
to establish  that more advanced versions
of 
$S$ can verify 
\eq{tangref}
for the cases where 1) $D$ denotes 
the more elaborate
\txl{1} 
{\it (rather than semantic tableaux)}
deduction method,
$~2)~~S$
 can
possess,
for any fixed i, 
 a finite-sized Group-2 scheme that enables it to 
prove the $i-$th kernelized image of each 
$\Pi_1^*$
theorem of
$A$,
 and 3)  $S$ can verify 
its own Level-1 consistency.

Such an argument
will corroborate
\thx{ttt5}'s
 claims 
about 
\eq{nrxefl} and
\eq{rdilute}'s
reflection principles.
(For instance, it applies to
\eq{nrxefl}
because $\Phi \,= \,\Phi^\oplus$
 for all
$\Delta_0^*$ or $\Sigma_1^*$ sentences.)

We will not provide more details concerning 
\thx{ttt5}'s proof here because the overall structure
of its justification is analogous to \cite{ww1}'s
corroboration of its ``TangRoot'' reflection
principle.
It is desirable to
 keep this technical report's
mathematical discussion as brief as possible,
so that our discourse can
focus, instead, on 
an 
 epistemological interpretation 
of the significance of self justification.


\begin{thebibliography}{99}

\small

\parskip 0pt



\baselineskip = 0.98\normalbaselineskip
\baselineskip = 1.0\normalbaselineskip





\bibitem{Ad2}
Z. Adamowicz,
``Herbrand Consistency and Bounded
Arithmetic'', 
{\it Fundamenta Mathematica}
171 (2002) pp. 279-292.



\bibitem{AZ1}
Z. Adamowicz and P. Zbierski,
``On Herbrand consistency in weak theories'', 
{\it Archive for Mathematical Logic}
40 (2001) pp. 399-413.



\bibitem{BS76}
A. Bezboruah and J. 
 Shepherdson,
``G\"{o}del's Second Incompleteness Theorem for Q'',
{\it Journal of Symb Logic} 41 (1976)  503-512.


\bibitem{Br94}
 S. Bringsjord, private communications during 1994 
suggesting that I publish my results about self-justifying
logics in two stages, with its first phase presenting
a series of mathematical results and second
stage
providing a  philosophical interpretation of such results.
Selmer thought this course was preferable because philosophical
interpretations always lend themselves to
wider levels of debates.
Thus following
Bringsjord's suggestion,
we had postponed discussing the philosophical interpertation of
\cite{ww93}-\cite{ww11}'s
 mathematical results
until the current report.

\bibitem{Bu86}
S. R. Buss, {\it Bounded Arithmetic,} (Ph D Thesis)
Proof Theory 
Notes \#3,  Bibliopolis 1986.



\bibitem{BI95}
S. R.  Buss and A. Ignjatovic,
``Unprovability of Consistency Statements in Fragments of
Bounded Arithmetic'', 
{\it Annals of Pure and Applied Logic} 74 (1995) pp. 221-244.


\bibitem{Da97}
John W. Dawson Jr.
{\it Logical Dilemmas: The Life and Work of
Kurt G\"{o}del,} A.K. Peter Press, 1997.

\bibitem{End}
H. B. Enderton,
{\it A Mathematical Introduction to Logic},
Acadmic Press, 2001.


\bibitem{Fi90}
M. Fitting, {\it  First Order Logic 
and
Automated Theorem Proving,}
SpringerVerlag 1996.



\bibitem{FW90}
M. L. Fredman and D. Willard, Blasting Through the
 Information Theoretic Barrier with Fusion Trees,
in
{\it  Proceedings of 22nd ACM
 Symposium on the
Theory of Computing (1990)} pp. 1--7 
(whose abridged proof-sketches were later published
by \cite{FW93} in a more detailed form).
 

\bibitem{FW93}
M. L. Fredman and D. Willard,
Surpassing the Information Theoretic Barrier with Fusion Trees,
{\it Journal of Computer and Systems Sciences,} 
47 (1993) pp. 424--436.


\bibitem{Fr79b}
H. M. 
Friedman, ``Translatability and Relative Consistency'',
 Ohio State  Tech Report, 1979.



\bibitem{Go31}
K. G\"{o}del,
`` \"{U}ber formal unentscheidbare S\"{a}tse der Principia
Mathematica und Verwandte Systeme I'',
{\it Monatshefte f\"{u}r Math. Phys.} 37 (1931) pp. 349-360.

\bibitem{Go5}
R. Goldstein, {\it Incompleteness The Proof and Paradox of
Kurt G\"{o}del}, Atlas Books of W. Norton and Company, 2005.

\bibitem{HP91}
P. H\'{a}jek and P. Pudl\'{a}k,
{\it Metamathematics of First Order Arithmetic,}
Springer Verlag 1991. 


\bibitem{Hil25}
D. Hilbert, Uber das Unendliche, Springer 1925.
(English Translation (1967) by J V Heijenoort in {\it From
Frege to G\"{o}del: A Source Book in 
Math
Logic,}
Harvard  Press, page 375)


\bibitem{HB39}
D. Hilbert and P. Bernays,
{\it Grundlagen der Mathematic}, Springer 1939.

\bibitem{Je71}
R. 
 Jeroslow, ``Consistency Statements in Formal Mathematics'', 
{\it Fund Math}  72(1971)  17-40.



\bibitem{Kl38}
S. C.
Kleene,
``On the Notation of Ordinal Numbers'',
{\it Journal  Symb Logic}
3 (1938), 
150-156.


\bibitem{Ko5}
L. A. Ko{\l}odziejczyk, 
Private Email Communications during November
of 2005.



\bibitem{KT74}
G. Kreisel and G. Takeuti,
``Formally Self-Referential Propositions  for Cut-Free Classical
Analysis'',
{\it Dissertationes Mathematicae}
118
(1974) pp. 1--55

\bibitem{Lo55} M. H. L\"{o}b, A Solution to a Problem by Leon Henkin,
{\it Jour  Symb Logic}
20(1955)  115-118



\bibitem{Mend}
E. Mendelson, {\it Introduction to Mathematical Logic},
CRC Press, 2010.




\bibitem{Ne86}
E. Nelson, {\it  Predicative Arithmetic,} Princeton Math Notes,
 1986.


    \bibitem{Pa71}
  R. Parikh, ``Existence and Feasibility in Arithmetic'',
{\it Journal   Symb Logic}
   36 (1971) 494-508



\bibitem{Pu85}
P. Pudl\'{a}k,
``Cuts, Consistency Statements\&Interpretations'',
{\it Journal Symb Logic}
50(1985)423-442

\bibitem{Ro67}
H. Rogers Jr. , {\it 
Recursive Functions and Effective 
Compatibility,} McGrawHill 1967.

 
 
\bibitem{So94}
R. M.  Solovay,  Private 
telephone  
conversations
in 1994
describing Solovay's generalization of one of  Pudl\'{a}k's  theorems 
\cite{Pu85},
using 
some methods
of Nelson and Wilkie-Paris \cite{Ne86,WP87}.
(The Appendix A of 
\cite{ww1} is a 
4-page approximate reconstruction of
this conversation.)

\bibitem{Stanford}
Stanford Encyclopedia Entry on ``Kurt G''{o}del,
See the second pargraph of its Section 2.2.4 for a
description of G\"{o}del's reaction to
Turing's paper \cite{Tu37}. 

\bibitem{Sv7}
V. \v{S}vejdar,
 ''An Interpretation of Robinson Arithmetic in its
Grzegorczjk's Weaker Variant''
{\it Fundamenta Informaticae} 81 (2007)
pp. 347-354.

\bibitem{Tu37}
A. M. Turing,
``On Computable Numbers with an Application to the
Entscheidungsproblem'',
Proceedings of the London Mathematical Soceity
(Series 2),  42 (1937) pp. 230-265.


\bibitem{Vi5}
A. Visser,
 ``Faith and Falsity'',
 {\it Annals of Pure and Applied Logic}
131 (2005) pp. 103-131.



\bibitem{WP87}
A. J. Wilkie and J. B. Paris,    
``On the Scheme of Induction for Bounded
Arithmetic'', {\it
Annals of Pure and Applied Logic} (35) 1987, 261-302

\vspace*{+0.2 em}

\bibitem{ww93}
D. Willard,
``Self-Verifying Axiom Systems'', {\it
Proceedings of the
   Third Kurt G\"{o}del
Colloquium}
(1993),
Springer-Verlag LNCS\#713,  pp. 325-336.



\bibitem{sp0}
---, 
``The Semantic Tableaux Version of the Second
Incompleteness Theorem Extends 
to 
Robinson's 
Arithmetic 
Q'', {\it Proceedings of   Tableaux 2000 Conf},
SpringerVerlag LNAI\#1847, 415-430 

\bibitem{ww1}
---, ``Self-Verifying  Systems, the Incompleteness
Theorem and the
Tangibility 
Principle'', in
{\it Journal of Symbolic Logic}
$~66~ (2001)\,$ pp. 536-596.


\bibitem{ww2}
---, 
``How to Extend The Semantic Tableaux And
Cut-Free Versions of the Second
Incompleteness Theorem 
Almost 
to 
Robinson's Arithmetic Q'', in $~$
{\it Journal of Symbolic Logic}
$~\,67~ (2002)~$ pp. 465--496.


 \bibitem{sp2}
 ---, 
 ``Some Exceptions for
 the 
 Semantic 
 Tableaux Version of the
 Second Incompleteness Theorem'', 
 {\it Proceedings of   Tableaux 2002 Conf},
 SpringerVerlag LNAI\#2381, 
 pp. 281--297.


\bibitem{wwlogos}
---, 
``A Version of the
Second Incompleteness Theorem For Axiom
Systems that Recognize Addition 
But Not Multiplication as a Total Function'',
{\it First Order Logic Revisited,
(Year 2003 Proceedings 
FOL-75 Conference),}
Logos Verlag (Berlin) 2004, pp. 337--368.




\bibitem{ww5}
---, 
``An Exploration of the Partial Respects in which an Axiom
System Recognizing Solely Addition as a Total Function Can
Verify Its Own Consistency'', 
{\it Journal of Symbolic Logic} 70  (2005) pp. 1171-1209. 

\bibitem{ww6}
---, 
``On the Available Partial Respects in which
 an Axiomatization
for Real Valued  Arithmetic Can  Recognize its 
Consistency'', 
{\it Journal of  Symbolic Logic} 71 (2006)
pp. 1189-1199.

\bibitem{wwapal}
---,  
``A Generalization of the Second Incompleteness 
Theorem and Some Exceptions to It''.
{\it Annals of Pure and Applied Logic}
141 (2006)
pp. 472-496.



\bibitem{ww7}
---,  
``Passive induction and a solution to a Paris-Wilkie 
question'',
{\it Annals of Pure and Applied Logic}
146(2007) 
pp. 124-149.

\bibitem{ww9}
---,  
``Some 
Special
Axiomizations for
I$\Sigma_0$ 
Manage to
Evade
the Herbrandized Version of the Second Incompleteness Theorem'',
{\it Information and Computation}
207(2009)  1078-1093


\bibitem{ww11}
--- ,  
``A Detailed Examination
of Methods for Unifying, Simplifying and Extending 
Several
Results About Self-Justifying Logics'',$~$
{\rm http://arxiv.org/abs/1108.6330}.



\bibitem{Yo5}
P. Yourgrau, {\it A World Without Time: The Forgotten Legacy of
G\"{o}del and Einstein}, Basic Books (2005).
(Our particular quotations from \cite{Yo5},
appearing in \textsection   \ref{sss1},
can be found on
its page 58.)


\end{thebibliography}
\end{document}